\documentclass[a4paper,11pt]{amsart}

\usepackage[utf8]{inputenc}
\usepackage[english]{babel}
\usepackage{hyperref}
\usepackage{amsmath,amsthm,enumerate}
\usepackage{amssymb}
\usepackage{tikz}
\usepackage{enumitem}

\usepackage{listings} % Include the listings package
\usepackage{xcolor} % Optional for coloring code

\usepackage{mathtools}
\usetikzlibrary{calc,decorations.pathmorphing,decorations.text,fixedpointarithmetic}
\usepackage{subfigure}
\usepackage{refcount}
\usepackage[hmargin=2.5cm,vmargin=3cm]{geometry}
\usepackage{graphicx,color}    % for latex
\usepackage{epsfig}
\usepackage{longtable}
\usepackage{caption}
\usepackage{comment}

\usepackage{tikz,tkz-graph}
\usetikzlibrary{calc}
\usetikzlibrary{backgrounds, shapes.geometric}
\usepackage{rotating}
\usetikzlibrary{positioning}
\usepackage{pgffor}

\newtheorem{theorem}{Theorem}[section]

\newtheorem{algorithm}[theorem]{Algorithm}

\newtheorem{proposition}[theorem]{Proposition}
\newtheorem{lemma}[theorem]{Lemma}
\newtheorem{setup}[theorem]{Setup}

\newtheorem{corollary}[theorem]{Corollary}

\newtheorem{definition}[theorem]{Definition}

\newcommand{\G}{\mathcal{G}}
\newcommand{\T}{\mathcal{T}}
\newcommand{\V}{\mathcal{V}}

\begin{document}

\title{Partial independent transversals in multipartite graphs}
\author{Penny Haxell$^1$}
\address{$^1$Department of Combinatorics and Optimization, University of Waterloo, Waterloo ON Canada N2L 3G1.}
\email{pehaxell@uwaterloo.ca}
\thanks{$^1$Supported in part by an NSERC Discovery Grant, and in part by National Science
Foundation Grant No. DMS-1928930, while the author was in
residence at the Simons Laufer Mathematical Sciences Institute in
Berkeley, California, during the Spring 2025 semester.}
\author{Arpit Mittal$^2$}
\address{$^2$Forsyth Central High School, Cumming, GA 30040}
\email{arpit.mittal.2.71@gmail.com}
\thanks{$^2$This research was conducted as a \textit{Research Experience for High School Students}, in which the second author participated in 2023--2025 under the guidance of the third author.}
\author{Yi Zhao$^3$}
\address{$^3$Department of Mathematics and Statistics, Georgia State University, Atlanta, GA 30303}
\email{yzhao6@gsu.edu}
\thanks{$^3$Partially supported by NSF grant DMS 2300346 and Simons Collaboration Grant 710094.}

%\author[1]{Penny Haxell\thanks{Supported in part by an NSERC Discovery Grant, and in part by National Science
%$Foundation Grant No. DMS-1928930, while the author was in
%residence at the Simons Laufer Mathematical Sciences Institute in
%Berkeley, California, during the Spring 2025 semester.}} 
%\author[2]{Arpit Mittal\thanks{This research was conducted as a %\textit{Research Experience for High School Students}, in which the second %author participated in 2023-2024 under the guidance of the third author.}}
%\author[3]{Yi Zhao\thanks{Partially supported by NSF grant DMS 2300346 and Simons Collaboration Grant 710094.}}

	%\affil[1]{\small Department of Combinatorics and Optimization,
         % University of Waterloo, Waterloo ON Canada N2L 3G1. E-mail:
         % pehaxell@uwaterloo.ca}
	%\affil[2]{\small Forsyth Central High School, Cumming, GA 30040}
    %    \affil[3]{\small Department of Mathematics and Statistics, Georgia State University, Atlanta, GA 30303}
\date{\today}
%\date{\bf {Jan 31, 2024}}

\begin{abstract} 
Given integers $r>d\ge 0$ and an $r$-partite graph, an independent $(r-d)$-transversal or $(r-d)$-IT is an independent set of size $r-d$ that intersects each part in at most one vertex.
We show that every $r$-partite graph with maximum degree $\Delta$ and parts of size $n$ contains an $(r-d)$-IT if $n> 2\Delta (1-\frac{1}{q})$, provided $q= \lfloor \frac{r}{d+1}\rfloor\ge \frac{4r}{4d+5}$. This is tight when $q$ is even and extends a classical result of Haxell in the $d=0$ case. When $q= \lfloor \frac{r}{d+1} \rfloor\ge \frac{6r+6d+7}{6d+7}$ is odd, we show that $n> 2\Delta(1-\frac{1}{q-1})$ guarantees an $(r-d)$-IT in any $r$-partite graph.
This is also tight and extends a result of Haxell and Szab\'o in the $d=0$ case. 
In addition, we show that $n> 5\Delta/4$ guarantees a $5$-IT in any $6$-partite graph and this bound is tight, answering a question of Lo, Treglown and Zhao.
\end{abstract}

\maketitle

\section{Introduction}\label{intro}
Let $G=(V,E)$ be a graph with partition $V=V_1 \cup \cdots \cup V_r.$ 
%Let $\Delta$ be the maximum degree of $G.$  
%\begin{definition}
An \emph{independent transversal} of $G$ is an independent set with exactly one vertex in each $V_i$. Given $0\le d< r$, an \emph{independent $(r-d)$-transversal} or \emph{$(r-d)$-IT} is an independent set with one vertex in each of $r-d$ parts of $G$.
We sometimes call an $(r-d)$-IT with $d>0$ a \emph{partial IT} and an independent transversal a \emph{full IT}. 
Independent transversals have found applications in many areas, including combinatorics (e.g. \cite{ADOV03,A18,CR15,LSM19,R11}), groups and rings (e.g. \cite{BEGHM08,C12}) and combinatorial optimization (e.g. \cite{AFS12,MR4141349}).
%linear arboricity \cite{A98}, strong chromatic numbers \cite{H04}, list coloring \cite{H01}, hypergraph matching \cite{AH00}, etc.

The existence of full independent transversals under maximum degree conditions 
was first studied by Bollob\'as, Erd\H{o}s, and Szemer\'edi \cite{BES} in 1975 in the complementary form. This problem was extensively studied and eventually solved by Haxell \cite{HA01}, Haxell and Szab\'o \cite{HS06}, and Szab\'o and Tardos \cite{ST06}, see \eqref{eq:nr}.
%Let $f(n, r)$ be the largest minimum degree $\delta(G)$ among all $r$-partite graphs $G$ with parts of size $n$ and which do not contain a complete subgraph on $r$ vertices. The authors of \cite{BES}

Bollob\'as, Erd\H{o}s, and Szemer\'edi also considered the existence of partial independent transversals (again in the complementary form). They \cite[Theorem 3.1]{BES} determined the maximum number of edges in $r$-partite graphs with $n$ vertices in each part and without a copy of $K_t$ (complete graph on $t$ vertices) for any $t\le r$. The corresponding problem under a minimum degree condition was studied by Lo, Treglown, and Zhao \cite{LTZ22}.

Given integers $n$ and $2 \le t \le r$, let $f(n,r,t)$ denote the largest minimum degree $\delta(G)$ among all $r$-partite graphs $G$ with parts of size $n$ and without a copy of $K_{t}$. A result of Bollob\'as, Erd\H{o}s, and Straus \cite{BESt} implies that $f(n, r, 3)= \lfloor r/2 \rfloor n$ for all $r\ge 3$ while \cite[Theorem 3.1]{BES} implies that $f(n,r,t+1) = \left(r- r/t\right)n$ whenever $t$ divides $r$. The authors of \cite{LTZ22} observed that Tur\'an's theorem implies that %for all $n,r,t$.
\begin{equation} \label{trivialbound}
\left(r-\left\lceil \frac{r}{t}\right\rceil \right)n \le f(n,r,t+1) \leq\left(r-\frac{r}{t}\right)n.
\end{equation}
The main result of \cite{LTZ22} determined $f(n,r,t+1)$ when $t$ divides $r+1$,  when $r\ge (3t-1)(t-1)$, and $f(n, r, 4)$ for all $r\ne 7$. %$f(n, 7, 5)= 5n$. %$r\equiv -1 \pmod{t-1}$ and

The smallest unknown case (with respect to $r$) is $f(n, 6, 5)$. In this paper we resolve this and many other cases of $f(n,r,t+1)$ when $t>r/2$ by studying the complementary problem
of finding partial independent transversals under maximum degree conditions. 
%Let $\Delta(n,r,t)$ be the smallest maximum degree among all $r$-partite graphs with part size $n$ and without $t$-IT. It is clear that
%\begin{equation}\label{delta and f}
%\Delta(n,r,t)+f(n,r,t)=(r-1)n.    
%\end{equation}
%In this paper we find it convenient to define another extremal function. 
Instead of $f(n, r, t)$, we find it convenient to use another extremal function.
For $r,D,\Delta \in \mathbb{N}$, let $n(r, D, \Delta)$ be the largest $n \in \mathbb{N}$ such that there is an $r$-partite graph $G$ with at least $n$ vertices in each part, maximum degree $\Delta$, and without an $(r-D+1)$-IT.

The two functions $f(n, r, t)$ and $n(r, D, \Delta)$ are closely related, as shown in Lemma \ref{f and n}. Under the $n(r,D,\Delta)$ notation, the results of \cite{HA01,HS06,ST06} are as follows: given integers $r\ge 2$ and $\Delta>0$, we have
\begin{align}
    \label{eq:nr}
    n(r, 1, \Delta)=\begin{cases}
			\left\lfloor 2\Delta \left(1- \frac{1}{r}\right) \right\rfloor, & \text{if $r$ is even}\\
             \left\lfloor 2\Delta \left(1- \frac{1}{r-1}\right) \right\rfloor, & \text{if $r$ is odd.}
		 \end{cases}
\end{align}

Given an $r$-partite graph $G$, let $G'$ be obtained from $G$ by adding $d$ vertex-disjoint copies of $K_r$ with one vertex in each part. Then $G$ has an $(r-d)$-IT if and only if $G'$ has a full IT. Using this standard idea for proving ``defect'' versions, one can adapt the proof of $n(r,1,\Delta) \le 2\Delta(1-1/r)$ (e.g. in \cite[Section 2]{H11}) to obtain that 
\begin{align} \label{eq:nrd1}
    n(r,d+1,\Delta)\leq 2\Delta\left(1-\frac{d+1}{r}\right).
\end{align}
Using multiple copies of \cite[Construction 3.3]{ST06} of Szab\'o and Tardos, one can easily see that \eqref{eq:nrd1} is tight when $r/(d+1)$ is an even integer. This already improves \eqref{trivialbound} and extends the first case of \eqref{eq:nr}. 

Our first result improves \eqref{eq:nrd1} substantially.
\begin{theorem}\label{genbound}
 Given integers $r, d\ge 0$, let $q,k$ be integers such that $r=q(d+1)+k$ and $0\le k\leq d$. If $G$ is an $r$-partite graph with maximum degree at most $\Delta$ and with vertex classes $V_1, \dots, V_r$ of size
  \[|V_i| >\max\left\{2\Delta\left(1-\frac{4d+5}{4r}\right),2\Delta\left(1-\frac{1}{q}\right)\right\},\] 
then $G$ has an $(r-d)$-IT.  
In other words, 
\[n(r,d+1,\Delta)\leq\max\left\{2\Delta\left(1-\frac{4d+5}{4r}\right),2\Delta\left(1-\frac{1}{q}\right)\right\}.\]
%For $r,d \in \mathbb{N}$ with $d \leqr$, we have 
%\[n(r,d+1)\leq\max\left\{2\left(1-\frac{4d+5}{4r}\right),2\left(1-\frac{1}{q}\right)\right\},\]
%where $q=\lfloor r/(d+1) \rfloor$.  
\end{theorem}
It is easy to see that the maximum in Theorem~\ref{genbound} is always less than or equal to $2\Delta\left(1-\frac{d+1}{r}\right)$. Thus, Theorem~\ref{genbound} is an improvement of \eqref{eq:nrd1}. Furthermore, together with the constructions in Section 2, Theorem~\ref{genbound} 
gives the value of $ n(r,d+1,\Delta)$ whenever $q\ge 4k$ is even (note that the maximum in Theorem~\ref{genbound} is equal to $2\Delta(1-1/q)$ if and only if $q\ge 4k$ or equivalently $q\ge \frac{4r}{4d+5}$). 

\begin{corollary}\label{smallk}
 Given integers $r, d\ge 0$, let $q,k$ be integers such that $r=q(d+1)+k$ and $0\le k\leq d$. If $q$ is even and $q\ge 4k$, then $n(r,d+1, \Delta)=\left\lfloor2\Delta\left(1-\frac{1}{q}\right)\right\rfloor$, and equivalently, $f(n, r, r-d)= (r-1)n - \left\lceil \frac{qn}{2(q-1)} \right\rceil$.
\end{corollary}

The $d=0$ case of Corollary \ref{smallk} says that $n(r, 1, \Delta)=\lfloor2\Delta(1-1/r)\rfloor$ for all even $r$, recovering the first case of \eqref{eq:nr} given by \cite{HA01,ST06}. 
%In addition, when $r=2t-1$ and $d=t-2$, 
%$q$ is even, $kq<t+1$, and $q \ge 4k$. 
%Indeed, we have $r=\frac{qt-k}{q-1}=q\left(\frac{qt-k}{q-1}-t\right)+k$. which implies $n(r,r-t)=2-2/q$, from which \eqref{f(r,t+1)} follows by applying Lemma \ref{f and n} and plugging in our expression for $r$.

\medskip
When $q=\lfloor \frac{r}{d+1} \rfloor$ is odd, we have the following results.

\begin{theorem}\label{oddq} 
Let $r,d, q \geq 0$ be integers such that $r=q(d+1)$. If $q\ge 3$ is odd, then
  \[n(r,d+1,\Delta)\leq\max\left\{ 2\Delta \left(1 -\frac{4d+5}{4r}\right),
  2\Delta \left(1-\frac{q}{q^2-1}\right)\right\}.
  \]
\end{theorem}

\begin{theorem}\label{odd q}
Given integers $r, d\ge 0$, let $q,k$ be integers such that $r=q(d+1)+k$ and $0\le k \leq d$. If $q \geq 3$ is odd, then
\[
n(r,d+1,\Delta)\leq\max\left\{2\Delta\left(1-\frac{6d+7}{6r}\right), 
2\Delta\left(1-\frac{1}{q-1}\right)\right\}.
\]
\end{theorem}

%Theorem~\ref{oddq} implies that $n(6,2) \leq 5/4$ and thus $f(6,5) \le 21/5$ by Lemma \ref{f and n}. Together with Corollary \ref{main construction} in Section~3, this gives the value of $f(6, 5)$.
Together with constructions in Section~2,  Theorem~\ref{oddq} determines $f(n, 6, 5)$, answering a specific question asked in \cite{LTZ22}. Previously, the best bounds were
\[
\frac{25}{6}n \le f(n, 6, 5)\le \frac{9}{2}n, 
\]
where the lower bound was given by \cite[Proposition 4.1]{LTZ22} and the upper bound follows from \eqref{trivialbound}. Theorem~\ref{genbound} or \eqref{eq:nrd1} gives $n(6, 2, \Delta)\le \frac43\Delta$, which implies $f(n, 6, 5)\le \frac{17}{4}n$ by Lemma \ref{f and n}.

\begin{corollary}\label{f(6,5)}
We have $n(6,2,\Delta)=\lfloor \frac54\Delta \rfloor$, and thus, $f(n, 6,5)=\lfloor \frac{21}{5}n \rfloor$.
\end{corollary}

Let us compare the bounds in Theorems~\ref{genbound}, \ref{oddq}, and \ref{odd q}.
Since %$1 -\frac{4d+5}{4r} \leq 1-\frac{6d+7}{6r}$ and 
$1-\frac{1}{q-1}< 1-\frac{q}{q^2-1}< 1-\frac{1}{q}$, the bound in Theorem~\ref{oddq} is at most the one in Theorem~\ref{genbound}. When $r=q(d+1)+k$ and $q\ge 6k$, the bound in Theorem~\ref{odd q} is at most the one in Theorem~\ref{genbound} because
\[
\max\left\{ 2 - \frac{6d+7}{3r}, 2 - \frac{2}{q-1}\right\} \leq 2 - \frac{2}{q} = \max\left\{2 - \frac{4d+5}{2r}, 2 - \frac{2}{q} \right\}.
\]
%Thus, whenever $q \geq 6k$, which is always true when Theorem \ref{oddq} applies, the upper bounds provided by Theorems \ref{oddq} and \ref{odd q} are always smaller than the bound in Theorem \ref{genbound}. %On the other hand, whether the bound in Theorem~\ref{oddq} or \ref{odd q} is smaller depends on the specific values for $r$ and $d$. In fact, 
Furthermore, when $k=0$ and $q\ge 6d+7 > \sqrt{4d+5}$, we have
\begin{align*}
    \max\left\{ 2 - \frac{6d+7}{3r}, 2 - \frac{2}{q-1}\right\} = 2-\frac{2}{q-1} < 2 - \frac{2q}{q^2-1} = \max\left\{ 2 -\frac{4d+5}{2r}, 2-\frac{2q}{q^2-1}\right\}.
\end{align*}
and thus, the bound in Theorem \ref{odd q} is smaller than the bound in Theorem \ref{oddq}. 

Together with the constructions in Section 2, Theorem \ref{odd q} determines the value of $ n(r,d+1,\Delta)$ whenever $q \ge 6d+6k+7$ (or equivalently $q\ge \frac{6r+6d+7}{6d+7}$) is odd.
\begin{corollary}\label{odd q tight}
Let $r, d\ge 0$ and $q,k$ be integers such that $r=q(d+1)+k$ and $0\le k \leq d$. If $q$ is odd and $q \geq 6d+6k+7$, then $n(r,d+1, \Delta)=\left\lfloor2\Delta\left(1-\frac{1}{q-1}\right)\right\rfloor$, and 
equivalently, $f(n, r, r-d)= (r-1)n - \left\lceil \frac{(q-1)n}{2(q-2)} \right\rceil$.
\end{corollary}

When $d=0$ (and thus $k=0$ and $q=r$), Corollary \ref{odd q tight} recovers the second case of \eqref{eq:nr} derived from the main result of \cite{HS06}, which assumes that $r\ge 7$ is odd.

Recall that $f(n, r, t)$ was determined in \cite{LTZ22} for $r=\Omega(t^2)$. Corollaries \ref{smallk} and \ref{odd q tight} together determine $n(r, d+1, \Delta)$ and $f(n, r, r-d)$ for $r=\Omega(d^2)$.
Indeed, given integers $d \geq 0$ and $r \geq 12d^2+20d+7$, we can write $r=q(d+1)+k$ such that $0\le k\le d$ and $q\ge 6d+6k+7\ge 4k$. Then, by Corollaries \ref{smallk} and \ref{odd q tight},
\begin{align*}
    n(r, d+1, \Delta)=\begin{cases}
			\left\lfloor 2\Delta \left(1- \frac{1}{q}\right) \right\rfloor, & \text{if $q$ is even}\\
             \left\lfloor 2\Delta \left(1- \frac{1}{q-1}\right) \right\rfloor, & \text{if $q$ is odd}.
		 \end{cases}
\end{align*}
%It follows from Corollaries \ref{smallk} and \ref{odd q tight} that for any integers In other words, we determine $n(r,d+1,\Delta)$ for $r=\Omega(d^2)$, which compares to \cite{LTZ22} where $n(r,d+1,\Delta)$ was determined for $r=\Omega((r-d)^2)$.                                                            
%Theorem~\ref{oddq} and Corollary~\ref{main construction} determine $f(n, 6, 5)$.
%shows that is sometimes tight, such as when $r=6$ and $d=1$. We can translate Theorem \ref{oddq} to the complementary setting to hence obtain tight results for $f$.
%\begin{corollary}
%Let $r=\frac{qt}{q-1}$ where $q \ge 3$ is odd. Then,
%\[f(r,t+1)\le \max\left\{\frac{q(4t-3)t}{(q-1)(4t-1)}-1, \frac{qt}{q-1}-\frac{q^2-1}{2q^2-2q-2}-1\right\}.\]
%In particular, $f(6,5)=21/5$.
%\end{corollary}

\subsection{Notation}
Given a graph $G$, let $V(G)$ and $E(G)$ be the edge and vertex sets of $G$ respectively. For $Z \subseteq E(G)$, we write $V(Z)$ to be the set of vertices incident to $Z$.  Given $v \in V(G)$ and $I \subseteq V(G)$, let $N(v,I)=\{ x \in I: \{x,v\} \in E(G)\}$ be the neighborhood of $v$ in $I$. %(note that $v \not\in N(v,I)$). 
Suppose $G$ is $r$-partite with parts $V_1, \cdots, V_r$. 
Given $I \subseteq V(G)$, let $S(I)=\{V_i \colon I \cap V_i \neq \emptyset\}$. We also let $\mathcal{G}_I$ be the multigraph obtained from the induced subgraph $G[I]$ by contracting all the vertices of $V_i \cap I$ into a single vertex denoted by $V_i$ (thus, $V(\mathcal{G}_I)=S(I)$). For $Y \subseteq[r]:= \{1,2,\ldots,r\}$, we let $G_Y$ be $G[\cup_{i \in Y} V_i]$. %We write $\Delta(G)$, or $\Delta$ when $G$ is understood, to mean the maximum degree of $G$.

\subsection{Organization}
We begin by considering constructions in Section~\ref{sec:cons}, where we prove Proposition~\ref{big construction} and Corollary \ref{main construction} and derive Corollaries \ref{smallk}, \ref{f(6,5)}, and \ref{odd q tight}. We prove Theorem \ref{genbound} in Section~\ref{sec:IMC} by generalizing the theory of \emph{Induced Matching Configurations} (IMCs) introduced in \cite{HS06}. In Section~\ref{sec:odd q results} we 
prove similar structural results as in \cite{HS06} and derive Theorems \ref{oddq} and \ref{odd q}. We give concluding remarks in the last section.
%which guarantee an $(r-d)$-IT when $r/(d+1)$ and $\lfloor r/(d+1) \rfloor$ are odd integers, respectively.

\section{Constructions and proofs of corollaries}
\label{sec:cons}

We first show how the two extremal functions $f(n,r, t)$ and $n(r, D, \Delta)$ defined earlier are related. To do so, it will be convenient to define another extremal function $\Delta(n,r,t)$ to be the smallest $\Delta \in \mathbb{N}$ such that there is an $r$-partite graph with at least $n$ vertices in each part, maximum degree $\Delta$, and without a $t$-IT. 

\begin{lemma}\label{f and n}
Suppose $r, d, n, \Delta$ are integers such that $r\ge 2$, $0\le d< r$, and $n, \Delta \ge 1$.
Let $c\ge 1$ be a real number. Then the following are equivalent:

\begin{enumerate}[label=\textup{(\roman*)}]
    \item $n(r,d+1,\Delta)=\lfloor c\Delta \rfloor$,
    \item $\Delta(n,r,r-d)=\lceil n/c \rceil$,
    \item $f(n,r,r-d)= (r-1)n - \lceil n/c \rceil$.
\end{enumerate}
\end{lemma}

\begin{proof}
We fix integers $r,d$ with $r\ge 2$ and $0 \le d < r$ throughout the proof. We first observe that for $n,\Delta \in \mathbb{N}$, 
\begin{equation}\label{n and delta}
n> \lfloor c\Delta \rfloor \quad \text{if and only if} \quad \Delta < \lceil n/c \rceil.   
\end{equation} %\textbf{YZ, this looks correct though the proof is not fully correct. We should label it as an equation and mention the equation whenever applying it.}
Indeed, as $n$ is an integer, $n>\lfloor c\Delta \rfloor$ implies $n \ge \lfloor c\Delta \rfloor+1 > c\Delta$, so $\Delta < n/c \le \lceil n/c \rceil$.
%\textbf{YZ: this does not give $\Delta < \lceil n/c \rceil$.}
On the other hand, $\Delta < \lceil n/c \rceil$ gives $\Delta \le \lceil n/c \rceil -1 < n/c$, so $n>c\Delta \ge \lfloor c\Delta \rfloor$. %We also have that $n=\lfloor c\Delta \rfloor$ is equivalent to $\Delta=\lceil n/c \rceil$. We can see this by noting that $n=\lfloor c \Delta\rfloor$ is equivalent to $c\Delta-1 < n \le c\Delta$, which holds if and only if $n/c \le \Delta<n/c+1/c<n/c+1$.

(i) $\implies$ (ii): We assume that $n(r,d+1,\Delta)=\lfloor c\Delta \rfloor$. 
Suppose $n\in \mathbb{N}$ and $G$ is $r$-partite such that each class has exactly $n$ vertices and $\Delta(G)=\Delta < \lceil n/c \rceil$. Then $n>\lfloor c\Delta \rfloor$ by \eqref{n and delta}. Since $n(r,d+1,\Delta)=\lfloor c\Delta \rfloor$, it follows that $G$ has an $(r-d)$-IT. 
Thus, $\Delta(n,r,r-d) \ge \lceil n/c \rceil$. 
On the other hand, the assumption that $n(r,d+1,\Delta)=\lfloor c\Delta \rfloor$ also implies that, for every $\Delta \in \mathbb{N}$, there exists an $r$-partite graph $G$ with maximum degree $\Delta$ and at least $\lfloor c\Delta \rfloor$ vertices in each part and without an $(r-d)$-IT. Now, given $n \in \mathbb{N}$, let $\Delta = \lceil n/c \rceil$ and $G$ be an $r$-partite graph with maximum degree $\Delta = \lceil n/c \rceil$ and at least $\lfloor c \lceil n/c \rceil \rfloor$ vertices in each part, that has no $(r-d)$-IT. Note that $\lfloor c \lceil n/c \rceil \rfloor \ge n$ because $c \lceil n/c \rceil  \ge c (n/c) = n$. 
%\[\left\lfloor c \left\lceil \frac{n}{c} \right\rceil \right\rfloor>c\left\lceil \frac{n}{c} \right\rceil-1 \ge n-1,\]
The existence of $G$ shows that $\Delta(n,r,r-d) \le \lceil n/c \rceil$. 

(ii) $\implies$ (i): Assume that $\Delta(n,r,r-d)=\lceil n/c \rceil$. Suppose $\Delta \in \mathbb{N}$ and $G$ is an $r$-partite graph with maximum degree $\Delta$ and parts of size $n>\lfloor c\Delta \rfloor$. We know $\Delta<\lceil n/c \rceil$ from \eqref{n and delta}.
Since $\Delta(n,r,r-d)=\lceil n/c \rceil$, it follows that $G$ has an $(r-d)$-IT, establishing $n(r,d+1, \Delta) \le \lfloor c\Delta \rfloor$. On the other hand, the assumption $\Delta(n,r,r-d)=\lceil n/c \rceil$ implies that, for every $n$, there exists an $r$-partite graph $G$ with parts of size at least $n$, maximum degree $\Delta(G)=\lceil n/c \rceil$, with no $(r-d)$-IT. Now, given $\Delta \in \mathbb{N}$, let $n=\lfloor c\Delta \rfloor$ and $G$ be the $r$-partite graph with parts of size at least $n=\lfloor c\Delta \rfloor$ and maximum degree $\lceil \lfloor c\Delta \rfloor/c \rceil$ and without an $(r-d)$-IT. Since $c\ge 1$, we have
\[
\Delta-1 \le \frac{c\Delta-1}{c}<\frac{\lfloor c\Delta \rfloor}{c}<\frac{c\Delta}{c}= \Delta,
\]
%\Delta-1 \le \frac{c\Delta-1}{c}<\left\lceil\frac{\lfloor c\Delta \rfloor}{c}\right\rceil<\frac{\lfloor c\Delta\rfloor}{c}+1\le\Delta+1,
and consequently, $\lceil \lfloor c\Delta \rfloor/c \rceil=\Delta$. 
The existence of $G$ shows that $n(r,d+1,\Delta) \ge \lfloor c\Delta \rfloor$.

(ii) $\iff$ (iii): By considering the complements of graphs, we have
\[\Delta(n,r,r-d)+f(n,r,r-d)=(r-1)n.\]
Hence, $\Delta(n,r,r-d)=\lceil n/c \rceil$ if and only if $f(n,r,r-d)= (r-1)n - \lceil n/c \rceil$.
\end{proof}

Since our results in this paper assume $d+1\le r/2$, the assumption $c\ge 1$ in Lemma~\ref{f and n} always holds. Indeed, given $r, \Delta\in \mathbb{N}$, 
let $G$ be the union of $\lfloor r/2 \rfloor$ vertex-disjoint copies of $K_{\Delta, \Delta}$ and an isolated set of $\Delta$ vertices if $r$ is odd. Then $\Delta(G)=\Delta$ and a maximum IT of $G$ has size $\lceil r/2 \rceil \le r - (d+1)$ if $d+1\le r/2$. This implies that $n(r, d+1, \Delta)\ge \Delta$.

\medskip
We now give several properties of $n(r, D, \Delta)$.
%For notational convenience, let $D$ be equal to $d+1$. 
It is clear that
%By dividing by $\Delta$ and taking the limit $\Delta\to\infty$, we obtain
\begin{align}
\label{eq:nrD}
n(r',D,\Delta)\le n(r,D,\Delta) \le n(r,D',\Delta) \quad \text{if} \quad r\ge r' \quad \text{and} \quad D\ge D'. 
\end{align}

By building new constructions from old ones, we derive the following proposition.
\begin{proposition}\label{big construction} The following holds for any positive integers $m, r, D, j, l$.
    \begin{enumerate}[label=\textup{(\roman*)}]
        \item\label{construction i} $n(r, D-1, \Delta)\ge n(r, D, \Delta) + \left\lfloor \frac{\Delta}{r-1}\right\rfloor$;
        \item\label{construction ii} $n(mr, mD, \Delta)\ge n(r, D, \Delta)$;
        \item\label{construction iii} $n(m r, (m-j) D, \Delta)\ge n(r, D, \Delta)+ \left\lfloor \frac{\Delta}{(l-1)r}\right\rfloor$ if $m=jl$;
        \item\label{construction iv} $n(m r, (m-j-1) D, \Delta)\ge n(r, D, \Delta)+ \left\lfloor\frac{\Delta}{(l-1)r}\right\rfloor + \left\lfloor\frac{\Delta}{(m-1)r}\right\rfloor$ if $m=jl$.
    \end{enumerate}
\end{proposition}
\begin{proof}
We say that a graph is an \emph{$(r,D)$-construction} if it has $r$ classes and its maximum IT has size at most $r-D$. Suppose $G$ is an optimal $(r, D)$-construction with maximum degree $\Delta$ and $r$ classes of size $n = n(r, D,\Delta)$. We will construct a $(r', D)$-construction $G'$ with maximum degree $\Delta$ and classes of size $n$, where $r'=r$ for \ref{construction i}, $r'=mr$ for \ref{construction ii}-\ref{construction iv}, and $D'=D-1, mD, (m-j)D$, and $(m-j-1)D$ for \ref{construction i}-\ref{construction iv} respectively.  

\smallskip
We first prove \ref{construction i}. Let $G'$ be the union of $G$ with $K_{r}(\lfloor \frac{\Delta}{r-1} \rfloor)$ (a at most $\Delta$-regular $r$-partite graph) such that each part of $G'$ has size 
\[
n + \left\lfloor \frac{\Delta}{r-1} \right\rfloor= n(r, D, \Delta) + \left\lfloor \frac{\Delta}{r-1} \right \rfloor.
\]
An IT of $G'$ has at most one vertex from $K_{r}(\lfloor\frac{\Delta}{r-1}\rfloor)$ and thus at most $r-D+1$ vertices in total. Thus $G'$ is an $(r, D-1)$-construction giving \ref{construction i}.

\smallskip
For \ref{construction ii}, let $G'$ be the union of $m$ disjoint copies of $G$. Then each part of $G'$ has $n(r, D, \Delta)$ vertices, as desired.

\smallskip
To prove \ref{construction iii}, we first show the case when $j=1$:
\begin{equation}\label{construction p2}
n(m r, (m-1) D, \Delta)\ge n(r, D, \Delta)+ \left\lfloor \frac{\Delta}{(m-1)r} \right\rfloor.
\end{equation}
%For convenience, write $K_{m\times r_0} = K_m(r_0)$.
To see \eqref{construction p2}, we take a union of $m$ (disjoint) copies of $G$ and a copy of $K_{m}\left(r\left\lfloor\frac{\Delta}{(m-1)r}\right\rfloor\right)$ such that each copy of $G$ is attached to an independent set of $r\left\lfloor\frac{\Delta}{(m-1)r}\right\rfloor$ vertices that are evenly distributed to $r$ classes. 
The resulting graph $G'$ has
\[
n + \left\lfloor\frac{\Delta}{(m-1)r}\right\rfloor = n(r, D, \Delta) + \left\lfloor\frac{\Delta}{(m-1)r}\right\rfloor
\]
vertices in each class.
It remains to show that a maximum IT of $G'$ has size at most $mr-(m-1)D$.
We regard each copy of $G$ together with $r\left\lfloor\frac{\Delta}{(m-1)r}\right\rfloor$ added vertices a \emph{row} and call 
the vertices of $G$ \emph{large} and the other vertices \emph{small}.
%We call the classes of $G_0$ \emph{large} and the classes of $K_{m\times r_0}(\frac{\Delta}{(m-1)r_0})$ \emph{small}. 
By definition, an IT of $G'$ misses at least $D$ large classes from each row and intersects small classes from at most one row.
%copy of $K_{1\times r_0}(\frac{\Delta}{(m-1)r_0})$. 
Thus the IT misses at least $mD - D$ classes of $G'$, confirming \eqref{construction p2}.

We now derive \ref{construction iii} from \ref{construction ii} and \eqref{construction p2}. Since $m- j = j(l-1)$, we have 
$n(m r, (m-j) D, \Delta)\ge n(l r, (l-1) D, \Delta)$ by \ref{construction ii}. By \eqref{construction p2}, we have $n(l r, (l-1) D, \Delta)\ge n(r, D, \Delta) + \left\lfloor\frac{\Delta}{(l-1)r}\right\rfloor$, proving \ref{construction iii}.

\smallskip
To see \ref{construction iv}, let $G'$ be the (disjoint) union of the following three graphs:
\begin{itemize}
    \item $m$ copies of $G$, arranged in $m$ \emph{rows}, and we call their vertices \emph{large};
    \item $j$ copies of $K_{l}\left(r\left\lfloor\frac{\Delta}{(l-1)r}\right\rfloor\right)$ such that each row has $r\left\lfloor\frac{\Delta}{(l-1)r}\right\rfloor$ independent vertices evenly distributed into 
    $r$ classes, %of size $\frac{\Delta}{(l-1)r_0}$, whom we call 
    and we call their vertices \emph{medium}; 
    \item one copy of $K_{m}(r\left\lfloor\frac{\Delta}{(m-1)r}\right\rfloor)$ such that each row has $r\left\lfloor\frac{\Delta}{(m-1)r}\right\rfloor$ independent vertices evenly distributed into 
    $r$ classes, and we call their vertices \emph{small}.
\end{itemize}
The resulting graph $G'$ has
\[
n + \left\lfloor\frac{\Delta}{(l-1)r}\right\rfloor + \left\lfloor\frac{\Delta}{(m-1)r}\right\rfloor =n(r, D, \Delta)+ \left\lfloor\frac{\Delta}{(l-1)r}\right\rfloor + \left\lfloor\frac{\Delta}{(m-1)r}\right\rfloor.
\]
vertices in each class.
It remains to show that the largest IT of $G'$ missed at least $(m- j-1)D$ vertices. Each IT of $G'$ has small vertices from at most one row, medium vertices from at most $j$ (additional) rows, and thus must use large vertices in at least $m-j-1$ rows. Since an IT of $G'$ misses 
at least $D$ vertices in each row, it thus misses at least $(m-j-1)D$ vertices in total, proving \ref{construction iv}.  
\end{proof}

\begin{corollary}\label{main construction}
Given integers $r\ge 1$ and $d\ge 0$, let $q, i, k$ be integers such that 
$q\ge 2$ is even, $1\le i\le d+2$, $0\le k< d+i$,  $r=q(d+i)+k$, and either $i=1$ or $i-1$ divides $d+i$.
Then,  
\[n(r,d+1, \Delta) \ge \left\lfloor 2\Delta\left(1-\frac{1}{q}\right) \right\rfloor + \left\lfloor\frac{(i-1)\Delta}{(d+1)q}\right\rfloor.\]
\end{corollary}
\begin{proof}
%We see that $n(r,D) \ge n(r-k, D)=n(q(d+i),D)$ by adding $k$ classes of $n(r-k,D)$ isolated vertices to an optimal $(r-k,D)$ construction with maximum degree $\Delta$ and class sizes of $n(r-k,D)$ 
%Observe that when $i=1$ the bound in Corollary \ref{main construction} is . 
Suppose $i=1$. Then we have
\[
n(q(d+1)+k,d+1, \Delta) \overset{\eqref{eq:nrD}}{\ge} n(q(d+1),d+1, \Delta) \overset{\ref{big construction}~\ref{construction ii}}{\ge} n(q,1, \Delta)\overset{\eqref{eq:nr}}{=}\lfloor2\Delta(1-1/q)\rfloor.
\]

Now assume $i>1$ and $d+i = (i-1)l$ for some integer $l$. Then,
\begin{align*}
    n(r,d+1, \Delta) &\overset{\eqref{eq:nrD}}{\ge} n(q(d+i),d+1, \Delta) \overset{\ref{big construction}~\ref{construction iii}}{\ge} n(q,1, \Delta)+\left\lfloor\frac{\Delta}{(l-1)q}\right\rfloor\\
    &\overset{\eqref{eq:nr}}{=} \left\lfloor 2\Delta\left(1-\frac{1}{q}\right) \right\rfloor + \left\lfloor\frac{\Delta}{\left(\frac{d+i}{i-1}-1\right)q}\right\rfloor=\left\lfloor 2\Delta\left(1-\frac{1}{q}\right) \right\rfloor + \left\lfloor\frac{(i-1)\Delta}{(d+1)q}\right\rfloor. \qedhere
\end{align*}
%Then by \cite{ST06} and Proposition~\ref{big construction}~(ii) with $m=d+i, r_0=q, j=i-1, l=m/j$, and $D_0=1$, we have $n(r,D) \ge n(q(d+i),D) \ge n(q,1)+1/(q(l-1))=(2-2/q)+1/(q((d+i)/(i-1)-1))=\left(2-\frac{2d+3-i}{q(d+1)}\right)$.
\end{proof}

Corollaries~\ref{smallk}, \ref{f(6,5)}, and \ref{odd q tight}
follow from Proposition~\ref{big construction} and Corollary~\ref{main construction} easily (assuming Theorems~\ref{genbound}, \ref{oddq}, and \ref{odd q}).

%Corollary \ref{main construction} with $i=1$ also provides the necessary lower bound on $n(r,d+1)$ required to establish Corollary \ref{smallk}.

\begin{proof}[Proof of Corollary~\ref{smallk}] Recall that $2-\frac{4d+5}{2r}\leq
   2-\frac2q$ precisely when $4k\leq q$. Hence 
   \[
   n(r,d+1, \Delta)\leq \left\lfloor \max\left\{2\Delta\left(1-\frac{4d+5}{4r}\right),2\Delta\left(1-\frac{1}{q}\right)\right\} \right\rfloor=\left\lfloor2\Delta-\frac{2\Delta}{q} \right\rfloor
   \]
   by Theorem \ref{genbound}. On the other hand, applying Corollary \ref{main construction} with $i=1$ gives $n(r,d+1, \Delta) \ge \lfloor2\Delta-2\Delta/q\rfloor$. Thus, $n(r,d+1, \Delta) = \lfloor2\Delta (1 - 1/q)\rfloor$. 
   By Lemma \ref{f and n}, we have  
   $f(n, r, r-d)= (r-1)n - \left\lceil \frac{qn}{2(q-1)} \right\rceil$.
\end{proof}

\begin{proof}[Proof of Corollary~\ref{f(6,5)}]
We have $n(6,2, \Delta) \le \lfloor5\Delta/4\rfloor$ by Theorem \ref{oddq} with $r=6$ and $d=1$. On the other hand, we have 
\[
n(6,2, \Delta) \ge n(2,1, \Delta) +\lfloor\Delta/4\rfloor = \lfloor5\Delta/4\rfloor
\]
by Proposition~\ref{big construction}~\ref{construction iii} with $m=3$, $j=1$, and $l=3$
(this can also be derived from Corollary \ref{main construction} with $q=2$ and $i=2$). 
Hence, $n(6,2, \Delta)= \lfloor5\Delta/4\rfloor$. Consequently, $f(n, 6,5)= 5n - \lceil 4n/5 \rceil=\lfloor 21n/5 \rfloor$ by Lemma \ref{f and n}.
\end{proof}

\begin{proof}[Proof of Corollary \ref{odd q tight}]
Since $q \ge 6d+6k+7$ is odd, Theorem \ref{odd q} gives that
\[
n(r,d+1) \le \left\lfloor\max\left\{2\Delta\left(1-\frac{6d+7}{6r}\right), 
2\Delta\left(1-\frac{1}{q-1}\right)\right\}\right\rfloor= \left\lfloor2\Delta-\frac{2\Delta}{q-1}\right\rfloor.
\]
On the other hand, since $r=q(d+1)+k$, we have
\[n(r,d+1, \Delta)\overset{\eqref{eq:nrD}}{\ge} n(q(d+1),d+1, \Delta) \overset{\ref{big construction}~\ref{construction ii}}{\ge} n(q,1, \Delta) \overset{\eqref{eq:nr}}{=} \left\lfloor2\Delta-\frac{2\Delta}{q-1}\right\rfloor.\]
Thus, $n(r,d+1, \Delta)=\lfloor2\Delta-2\Delta/(q-1)\rfloor$.
By Lemma \ref{f and n}, it follows that  
   $f(n, r, r-d)= (r-1)n - \left\lceil \frac{(q-1)n}{2(q-2)} \right\rceil$, as desired.
\end{proof}

The proof of Corollary~\ref{f(6,5)} shows that Proposition~\ref{big construction}~\ref{construction iii} is tight. 
%Note that Corollary \ref{smallk} gives us examples of when Proposition~\ref{big construction}~\ref{construction ii} is tight. 
Proposition~\ref{big construction} and Corollary~\ref{smallk} together imply that, for any even $q\ge 2$ and $0\le k\le q/4$,
\begin{align*}
    \left\lfloor2\Delta - \frac{2\Delta}{q}\right\rfloor &\overset{\ref{smallk}}{=} n(q(d+1)+k, d+1, \Delta)\overset{\eqref{eq:nrD}}{\ge} n(q(d+1), d+1, \Delta)\\
    &\overset{\ref{big construction}~\ref{construction ii}}{\ge} n(q, 1, \Delta) \overset{\eqref{eq:nr}}{=} \left\lfloor2\Delta - \frac{2\Delta}{q}\right\rfloor.
\end{align*}
In particular, this shows that the function $n(r, D)$ is not strictly increasing in $r$ and Proposition~\ref{big construction}~\ref{construction ii} is tight. (Note that, in contrast, Proposition~\ref{big construction}~\ref{construction i} shows that $n(r, D)$ is strictly decreasing in $D$.)
%When $i=2$ we can translate Corollary \ref{main construction} as follows. Let $r=\frac{q(t-1)-k}{q-1}$ where $k \le (t-1)/q$ and $q \ge 2$ is even. Then
%\begin{equation}\label{f construction}
%f(r,t+1) \ge \frac{q(t-1)-k}{q-1}-\frac{q(t-q-k)}{2(q-1)(t-q-k)+(q-1)}-1.    
%\end{equation}

\section{Induced Matching Configurations and proof of Theorem~\ref{genbound}}
\label{sec:IMC}
Suppose $G$ is an $r$-partite graph with vertex classes $V_1, \dots, V_r$. Let $\V=\{ V_1, \dots, V_r\}$.
%Given a set of vertices $I \subseteq V(G)$, we say that $I=I_1 \cup \cdots \cup I_s$ is a partition of $I$ into $s$ components if $N(I_k) \cap I_j=\emptyset$ and $I_j \cap I_k=\emptyset$ for $i \neq j$. 
Fix a subset $I \subseteq V(G)$. Recall that $S(I)$ is the set of classes that intersect $I$, and $\G_I$ is the multigraph formed by contracting all vertices of $G[I]$ into one vertex, which we still denote by $V_i$. Note that $V(\G_I) = S(I)$ and $\G_I$ has parallel edges if $G[I]$ has multiple edges between two classes $V_i$ and $V_j$.
%$V(G_I)=\{v_i \colon V_i \cap I \neq \emptyset\}$ and $\{v_i,v_j\}\in E(G_I)$ iff there exist $w_i \in I \cap V_i, w_j \in I \cap V_j$ such that $\{w_i,w_j\} \in E(G)$. 
\begin{definition}
A set of vertices $I \subseteq V(G)$ is an Induced Matching Configuration (IMC) if $G[I]$ is a perfect matching and $\mathcal{G}_I$ is a forest. We say that $I$ is an IMC of $p$ components if $\G_I$ has $p$ components.
%$I$ admits a partition into $(d+1)$ components $I=I_1 \cup \cdots \cup I_{d+1}$ such that $G[I]$ is a perfect matching and for each component $I_k$, $G_{I_k}$ is a tree.
\end{definition}
If $I$ is an IMC, then every edge of $G[I]$ corresponds to a (unique) edge of $\G_I$ (because $\G_I$ is a forest and thus has no multiple edges). 
%$G[I]$ has at most one edge between any two classes  
Therefore, if $|I| = 2t$ and $|S(I)|=s$, then $\G_I$ has $s$ vertices and $t$ edges, thus consisting of $s-t$ components. We remark that IMCs were defined differently in \cite{HS06}, in which $\G_I$ is a tree on $\V$. 

%All the IMC's in 
%We remark that if $I$ contains vertices from $t$ distinct classes, then $|I|=2(t-d-1)$. In the following writing, when talking about a vertex $v_i \in V(I)$, we write $w_i$ to mean the (unique) neighbor of $v_i$ in $I$. We also sometimes omit the $(r-d)$ and refer to an $(r-d)$-IMC as an IMC.
%The following lemma is from \cite{HS06}, which was proved for IMCs. It allows us to extract a (partial) IT from a given IMC. 

The forest structure of $\G_I$ makes it easy to find a partial IT in $G[\cup_{V_i \in S(I)} V_i]$ covering all but one class of every component of $\G_I$.  
The following lemma is essentially \cite[Lemma 2.1]{HS06} but revised slightly due to our definition of IMCs. 
\begin{lemma}[{\hspace{1sp}\cite[Lemma 2.1]{HS06}}]\label{it from imc}
%Let $G$ be an $r$-partite graph with IMC $I$ with $|S(I)|=r$. Then, for any $i \in [r]$, $G$ has a partial IT $T$ of size $(r-1)$ contained in $I$ such that for $j \in [r]$, $|T \cap V_j|=1$ iff $j \neq i$.
Suppose $I$ is an IMC in $G$ and $I'$ is a subset of $I$ such that $\G_{I'}$ is a tree on $S\subseteq \V$. Then, for any $V_i\in S$, there is an $(|S|-1)$-IT on $S\setminus \{V_i\}$. Further, if $v\in V_i\in S$ is not dominated by $I'$, then there is an $|S|$-IT on $S$.   
\end{lemma}

\begin{comment}
\begin{proof}
Fix a component $C_k$ and let $c_k=C_k \cap C$. Then, considering $G[I]$ to be a rooted tree with $V(c_k)$ as the root, for all $j \in C_k-\{V(c_k)\}$ we add to $T$ the element in $V_j \cap I$ whose neighbor in $G_I$ is the parent class of $V_j$ in $G_I$. Thus, for every component $C_k$ we add $|C_k|-1$ vertices to $T$, so $|T|=t-s$. 
\end{proof}
\end{comment}
We need the concept of feasible pairs introduced in \cite{HS06} in order to prove Theorem~\ref{imc exists}. 
\begin{definition}\label{feasible pair}
Given an $r$-partite graph $G$, a pair $(I,T)$ with 
$I,T \subseteq V(G)$ is feasible if the following conditions hold:
\begin{enumerate}[label=(\alph*)]
    \item $T$ is a partial IT in $G$ of maximum size.
    \item $S(I \cap T)=S(I) \cap S(T)$, i.e. if $v \in T$ and $S(\{v\}) \in S(I)$, then $v \in I$.
    \item $G[I]$ is a forest, whose components are stars  with centers in $W:=I \setminus T$ and at least one leaf.
    %$|W|$  vertex disjoint components of the form $N(v,T)$ for $v \in W$,
    \item $\G_I$ is a forest on $S(I)$.
    \item Let $\T$ be the set of all partial independent transversals of $G$ on $S(T)$.
    %such that $T'\cap V_i\ne \emptyset$ if and only if $V_i\in S(T)$.
    Then, for every $v \in W,$ there is no $T' \in \T$ such that  $T' \cap W = \emptyset$, $|N(v,T')|<|N(v,T)|$, and $N(w,T')=N(w,T)$ for all $w \in W\setminus \{v\}$.
\end{enumerate}
\end{definition}

Note that if $T$ is a partial IT in $G$ of maximum size, then $(\emptyset, T)$ is a feasible pair.

The following algorithm allows us to construct feasible pairs that dominate all vertex classes intersecting them. Recall that $N(I)$ is the set of \emph{neighbors} of $I$, in other words, the vertices that are dominated by $I$.
%Given an $r$-partite graph $G$ on $V_1 \cup \cdots \cup V_r$ and $T\subseteq V(G)$, let $R= R(T) = \{V_1, \ldots, V_r\} \setminus S(T)$. 
% \begin{algorithm} 
% Start with a feasible pair $(I_0, T_0)$ in $G=V_1 \cup \cdots \cup V_r$. 
% %set $I\subset V(G)$ and an $(r-d-1)$-IT $T_0$ of $G$. 
% Initialize $I=I_0$ and $T=T_0$. Throughout the algorithm, we maintain $W=I\setminus T$ and $\T$ which is defined to be the set of all maximum partial ITs $T'$ on $S(T_0)$ such that $T'\cap W=\emptyset$, $N(v,T)=N(v,T_0)$ for all $v \in W$. 
% \begin{enumerate}
% \item[Step 1] Select a vertex $w\in \bigcup_{V_i\in R} V_i\setminus N(I)$ and $T'\in \T$ such that $\deg(w, T')$ is minimal. Update $I$ by adding $\{w\}\cup N(w, T')$ and update $T= T'$. Go to Step 2.
% \item[Step 2] If $I$ dominates all the vertices in $S(I)\cup R$, then stop and return $(I, T)$. %Update $S(I), W,$ and $\T$ accordingly. 

% If $I$ dominates $V_I= \bigcup_{V_i\in S(I)} V_i$ but not all the vertices in $R$, go to Step 1.
 
% If $V_I$ is not dominated by $I$, then select a vertex $w\in V_I\setminus N(I)$ and $T'\in \T$ such that $\deg(w, T')$ is minimal. Update $I$ by adding $\{w\}\cup N(w, T')$, update $T= T'$, and repeat Step 2.
% \end{enumerate}
% \end{algorithm}

% \textcolor{red}{The following is the revised version of the algorithm.}

\begin{algorithm}\label{alg} 
Start with a feasible pair $(I_0, T_0)$ in an $r$-partite graph $G$ on $V_1 \cup \cdots \cup V_r$. Let $R= \{V_1, \ldots, V_r\} \setminus S(T_0)$. 
%set $I\subset V(G)$ and an $(r-d-1)$-IT $T_0$ of $G$. 
Initialize $I=I_0$ and $T=T_0$. Throughout the algorithm, we maintain $W=I\setminus T$ and $\T$, the set of all maximum partial ITs $T'$ on $S(T)$ such that $T'\cap W=\emptyset$, $N(v,T')=N(v,T)$ for all $v \in W$. 
\begin{enumerate}
\item[Step 1] If $I$ dominates all the vertices in $S(I)\cup R$, then stop and return $(I, T)$
\item[Step 2] If $I$ dominates $V_I= \bigcup_{V_i\in S(I)} V_i$ but not all the vertices in $R$, go to Step 3.
 
If $V_I$ is not dominated by $I$, then select a vertex $w\in V_I\setminus N(I)$ and $T'\in \T$ such that $\deg(w, T')$ is minimal. Update $I$ by adding $\{w\}\cup N(w, T')$, update $T= T'$, and go to Step 1.
\item[Step 3] Select a vertex $w\in \bigcup_{V_i\in R} V_i\setminus N(I)$ and $T'\in \T$ such that $\deg(w, T')$ is minimal. Update $I$ by adding $\{w\}\cup N(w, T')$ and update $T= T'$. Go to Step 1.
\end{enumerate}
\end{algorithm}

\begin{lemma}\label{algorithm terminates}
For any input $(I_0, T_0)$, Algorithm \ref{alg} terminates and returns a feasible pair.
\end{lemma}

As the proof of Lemma \ref{algorithm terminates} is very similar to the corresponding part of the proof of \cite[Theorem 2.2]{HS06}, we defer it to Appendix \ref{appendix proof of lemma}. Indeed, 
Algorithm~\ref{alg} is the same as the one given in the proof of \cite[Theorem 2.2]{HS06}, except that in  \cite{HS06} the algorithm terminates when $I$ dominates $S(I)$, while Algorithm~\ref{alg} terminates when $I$ dominates $S(I)\cup R$. As a result, the resulting $\G_I$ is a forest, instead of a tree as in \cite{HS06}.

\begin{theorem}\label{imc exists}
Let $G$ be an $r$-partite graph on $\V= \{V_1, \dots, V_r\}$. Suppose that the largest partial IT of $G$ has size $r-d-1$. Given a feasible pair $(I_0, T_0)$, there exists a feasible pair $(I, T)$ in $G$ such that
\begin{enumerate}[label=\textup{(\roman*)}]
    \item $I_0\subseteq I$, $S(I)\ge 2$, $S(T)=S(T_0)$, and $T\cap V_i = T_0\cap V_i$ for every $V_i\in S(I_0)$.
    \item $I$ dominates all of the vertices in $S(I)\cup R$, where $R=\V\setminus S(T_0)$. 
    \item Let $F_I= (R\cup S(I), E(\G_I))$ be the extension of $\G_I$ to $R\cup S(I)$ if $R\not\subseteq S(I)$ (otherwise $F_I=\G_I$). Then $F_I$ is a forest of $d+1$ components, where each component contains exactly one $V_i\in R$ (which may be the only $V_i$ in that component).
    \item Let $t= |R\cup S(I)|$. Then $|I| \le 2(t -d -1)$.
\end{enumerate}
If $n > 2\Delta \left( 1 - \frac{2d + 3}{2r} \right)$, 
then $I$ is an IMC of $G$ and $|I| = 2(t -d -1)$.
\end{theorem}

\begin{comment}
    Let $G$ be an $r$-partite graph whose largest partial IT has size $r-d-1$. If $n > 2\Delta \left( 1 - \frac{2d + 3}{2r} \right)$, 
there exists an IMC $I\subset V(G)$ and an $(r-d-1)$-IT $T$ such that
\begin{enumerate}
\item Let $W= \{V_i: V_i\cap T = \emptyset\}$. Then $I$ dominates $W\cup S(I)$. 
\item $G[I]$ is an matching, and $e_G(V_i\cap I, V_j\cap I)\le 1$ for any $i\ne j\in [r]$.
\item Let $F_I= (W\cup S(I), E(\G_I))$ be the extension of $\G_I$ to $W\cup S(I)$ if $W\not\subseteq S(I)$ (otherwise $F_I=\G_I$). Then $F_I$ is a forest of $d+1$ components, where each component contains exactly one $V_i\in W$ (which may be the only vertex in the component).
\item Let $t= |W\cup S(I)|$. Then $|I| = 2(t -d -1)$.
\item $T\setminus I$ is a partial IT that is independent of $I$.
\end{enumerate}
Moreover, if $n>2\Delta \left( 1 - \frac{d + 1}{r - 1} \right)$ as well, then $t=r$, or equivalently, $|S(I)|=r-d-1$.
\end{comment}

\begin{proof}
Let $(I, T)$ be the feasible pair obtained by applying Algorithm \ref{alg} to $(I_0, T_0)$. We now show that $(I, T)$ satisfies the conditions in the statement of Theorem \ref{imc exists}.

Part (i) follows from the definition of feasible pairs and the algorithm. 

Part (ii) follows from the algorithm immediately.

To see Part (iii), we note that each tree of $\G_I$ has exactly one class from $R$, and $|R|=d+1$ because $|T_0|=|T|= r-d-1$. It is possible that a class $V_i\in R$ contains no vertex of $I$ (if all the vertices of $V_i$ are dominated by $I$ at some stage of the algorithm). 

For Part (iv), since $F_I$ is a forest on $t$ vertices with $d+1$ components, we have $|E(F)|= t - d -1$ and $|I|\le 2|E(F)|= 2(t-d-1)$.

We now show that if $n > 2\Delta \left( 1 - \frac{2d + 3}{2r} \right)$, then $G[I]$ is a perfect matching, or equivalently $|I|=2(t-d-1)$. Suppose that $|I|\le 2(t-d-1)-1$. Since $I$ dominates $S(I)\cup R$, we have $|I|\Delta\ge |S(I)\cup R| n$, which gives $\Delta(2(t-d-1)-1)> 2t\Delta\left(1-\frac{2d+3}{2r}\right)$. This implies $\Delta(2d+3)< \Delta\frac{t}{r}(2d+3)$, contradicting the fact that $r\ge t$.
%It remains to be shown that if $n>2\Delta\left(1-\frac{d+1}{r-1}\right)$, then $t=r$. Since $I$ dominates $S(I)\cup W$ and $|I|=2(t-d-1)$, we have  $2(t-d-1) \Delta > 2t\Delta\left(1-\frac{d+1}{r-1}\right)$. This implies that $t\frac{d+1}{r-1}> d+1$, or $t> r-1$, as desired. 
\end{proof}

Since we will make use of it later in this section, we explicitly state the $d=0$ case of Theorem \ref{imc exists}, which is precisely the content of \cite[Theorem 2.2]{HS06}

\begin{lemma}[{\hspace{1sp}\cite[Theorem 2.2]{HS06}}]\label{lem:noit}
Let $G$ be an $r$-partite graph on $V_1 \cup \ldots \cup V_r$.  Suppose that $G$ has no IT. Then, there exists a set of class indices $S \subseteq [r]$ and a nonempty set of edges $Z$ of the subgraph $G_S:= G[\bigcup_{i\in S} V_i]$ such that $V(Z)$ dominates $G_S$ and $|Z| \leq|S|-1$. Moreover, $Z \cap V_s \neq \emptyset$ for every $s \in S.$ 
\end{lemma}

When Algorithm~\ref{alg} terminates, the set $T\setminus I$ must be a partial IT that is independent of $I$ because if $w\in T\setminus I$ is adjacent to some $v\in I$, then $w$ would have been included in $I$ when $v$ was added to $I$.
The following proposition strengthens this fact by showing that such a partial IT exists outside $N_G(x)$ for any fixed $x\in V(G)$.

\begin{lemma}\label{prop b}
Let $G$ be an $r$-partite graph with parts $V_1, \dots, V_r$ of size $n>\Delta\left(2-\frac{2d+3}{r}\right)$. Let $I\subset V(G)$, $d+1\le t\le r$, and $U\subseteq \{V_1, \dots, V_r\}$ be such that $|U|= r -t$, $|I| \le 2(t-d-1)$, and $I$ dominates all the classes not in $U$.
%be a set of at most $2(t-d-1)$ vertices that dominates $t$ classes and let $U$ be the set of the classes not dominated by $I$. 
%$I$ be the IMC provided by Theorem~\ref{imc exists} and let $U=\{V_i: V_i\not\in W\cup S(I)\}$. 
For any vertex $x \in V(G)$, there exists an IT $T'$ of $G[U]$  such that $T' \cap N_G(I \cup \{x\}) = \emptyset$. 
\end{lemma}
\begin{proof}
Let $V'_i := V_i\setminus N_G(I \cup \{x\})$ for all $i$ and $U' = \{V'_i: V_i\in U\}$. Note that it is possible for some $V'_i=\emptyset$.
If there is an IT of $G[U']$, then we are done.
Suppose there is no IT of $G[U']$. If all $V'_i$ are nonempty, then by Lemma~\ref{lem:noit}, there exists a subset $U'' \subseteq U'$ such that $G[U'']$ is dominated by a set $X \subseteq \bigcup_{V'_i\in U''} V'_i$ of most $2(|U''|-1)$ vertices. If some $V'_{i_0}= \emptyset$, then let $U''=\{V_{i_0}\}$ and $X= \emptyset$.
In either case, the set $I \cup N_G(x) \cup X$ dominates $\bigcup \{V_i: V_i\not\in U$ or $V'_i\in U''\} $. Since $|I|\le 2(t-d-1)$, we have 
\[|I \cup \{x\} \cup X| \leq2(t-d-1)+2(|U''|-1)+1=2(t-d+|U''|-2)+1.\]
Since there are $(t+|U''|)n$ vertices in $\bigcup \{V_i: V_i\not\in U$ or $V'_i\in U''\}$, %$\bigcup_{V_i\not\in U, V'_i\in U''} V_i$, 
we derive that  
\[(2(t-d+|U''|-2)+1)\Delta \ge n(t+|U''|) > \Delta\left(2-\frac{2d+3}{r}\right)(t+|U''|),\] 
which implies $2(-d-2)+1>-\frac{2d+3}{r}(t+|U''|)$. This is a contradiction because $r = t + |U'| \ge t+|U''|$. 
\end{proof}

%Let $\mathcal{S}\subseteq \{V_1, \dots, V_r\}$ and 
Suppose $I\subseteq V(G)$ is a set that dominates all the classes of $R\cup S(I)$.
For every vertex $v \in I$, we let $A_v$ or $A_v(I)$ be the set of all vertices in the classes of $R\cup S(I)$ that are adjacent to $v$ but no other vertex in $I$. For example, if $v,w \in I$ are adjacent, then $w \in A_v$ and $v \in A_w$. By definition, if $x\in \bigcup \{V_i: V_i\in R\cup S(I) \}\setminus \bigcup_{v\in I} A_v$, then $x$ is adjacent to at least two vertices of $I$. In this next lemma, we bound the number of vertices in $G$ adjacent to at least two vertices of $I$ and the number of vertices in a given number of $A_v$'s using similar arguments as the ones for \cite[Lemma 3.1 (iii) and (iv)]{HS06}. 

\begin{lemma}\label{avsize}
Let $G$ be an $r$-partite graph with parts $V_1, \dots, V_r$ of size $n$. Suppose $I\subseteq V(G)$ is a set of size at most $2(t-d-1)$ that dominates all the classes of $R\cup S(I)$ with $t=|R\cup S(I)|$. Then
\begin{enumerate}[label=\textup{(\roman*)}]
    \item\label{odd q max D} The number of vertices $x\in V(G)$ such that $|N(x) \cap I| \geq 2$ is at most $2\Delta(t-d-1)-tn$,
    \item\label{odd q min Av} For any subset $Y\subseteq I$, we have $|\bigcup_{v\in Y} A_v| \ge (|Y|+4d+4 - 4t)\Delta + 2t n$,
    \item\label{old avsize} If $n>2\Delta\left(1-\frac{4d+5}{4r}\right)$, then for any $Y \subseteq I$, we have $|\cup_{v \in Y}A_v|>(|Y|-1)\Delta$.
\end{enumerate}
\end{lemma}

\begin{proof} Let $V'= \bigcup \{V_i: V_i\in S(I)\cup R\}$.

Part (i):  Let $D$ be the set of vertices $x\in V(G)$ such that $|N(x) \cap I| \geq 2$. Then $V'\setminus D= \bigcup_{v\in I} A_v$ consists of the vertices of $V'$ that are dominated by $I$ exactly once.
Hence 
\[
2\Delta(t-d-1)\ge \Delta|I| \ge 2|D|+(tn-|D|) \geq tn+|D|,
\]
which implies that $|D|\le 2\Delta(t-d-1)-tn$.

\smallskip
Part (ii): We know $\bigcup_{v\in Y} A_v = V'\setminus (D\cup \bigcup_{v\in I\setminus Y} A_v)$. Applying Part (i), we obtain that 
\begin{align*}
|\bigcup_{v\in Y} A_v| &\ge tn - (2\Delta(t-d-1)-tn) - (2(t-d-1)-|Y|) \Delta\\
&= (|Y|+4d+4 - 4t)\Delta + 2tn.
\end{align*}

\smallskip
Part (iii): By Part (ii), we have 

\begin{align*}
|\bigcup_{v\in Y} A_v| &> (|Y|+4d+4 - 4t)\Delta + 4t\Delta\left(1-\frac{4d+5}{4r}\right)\\
&=\left(|Y|+4d+4-\frac{t(4d+5)}{r}\right)\Delta
%\\&
\geq (|Y|-1)\Delta
\end{align*}
as $t\le r$.
\end{proof}

%Let $I$ be an IMC returned by Theorem~\ref{imc exists}. Recall that $F_I = (R\cup S(I), E(\G_I))$, as defined in Theorem~\ref{imc exists}~(iii) is a forest of $d+1$ components. where each component $C$ is a subset of $R\cup S(I)$ such that $F_I[C]$ is a tree. Note that it is possible that some component contains only one $V_i\in R$ such that $V_i\cap I= \emptyset$.

For convenience, we introduce the following setup for the next few lemmas. 
Part~\ref{setup}~(i) includes all the properties of an IMC returned by Theorem \ref{imc exists} while Part (ii) has a stronger lower bound for $n$ required by Lemma~\ref{avsize} \ref{old avsize} \ref{old avsize}.
%It is possible that $I$ dominates some $V_i\in R$ with $V_i\cap I= \emptyset$. After including such $V_i$ as a single-vertex component, 
%For simplicity, we will say that \emph{$I$ is an IMC of $d+1$ components} with components $I\cap C$ for each component $C$ of $F_I$.
%In particular,  

\begin{setup}\label{setup}
We have the following two setups:
\begin{enumerate}[label=\textup{(\roman*)}]
    \item Let $G$ be an $r$-partite graph on $V_1 \cup \cdots \cup V_r$. Let $T$ be a maximum IT of $G$ of size $r-d-1$. Assume that $\Delta(G)\le \Delta$ and $|V_i|\ge n>2\Delta\left(1-\frac{2d+3}{2r}\right)$ for all $i$. Let $I$ be an IMC of $G$ returned by Theorem \ref{imc exists}, which has the following properties:
    \begin{itemize}
        \item $I$ dominates all the vertices in $R \cup S(I)$, where $R= \{V_1, \dots, V_r\}\setminus S(T)$;
        \item $F_I= (R\cup S(I), E(\G_I))$ is a forest of $d+1$ components, where each component $C$ is a subset of $R\cup S(I)$ such that $F_I[C]$ is a tree.
        \item $|R \cup S(I)|=t$ and $|I|=2(t-d-1)$ for some $d+1\le t\le r$.
    \end{itemize} 
    \item Let $G$ be as in Setup~\ref{setup} (i) with the additional assumption $n>2\Delta\left(1-\frac{4d+5}{4r}\right)$.
    %$n>2\Delta-(4d+5)\Delta/2r$.
\end{enumerate}
\end{setup}

Suppose $G$ is as in Setup~\ref{setup} (i). If a vertex $x\in V_i$ and $V_i\in C$ for some component $C$ of $F_I$, then we may simply say $x$ is in $C$. In other words, we partition the vertices $G$ in $R\cup S(I)$ based on the components of $F_I$. In particular, 
%We say that \emph{two vertices of $G$ are in the same (or different) components} if the classes that contain them are in the same (or different) components of $F_I$.
we say that \emph{$I$ is an IMC of $d+1$ components}. For example, any two adjacent vertices $v, w\in I$ belong to the same component. 

In next few lemmas we may replace the IMC $I$ in Setup~\ref{setup} by another IMC $I'$ of $G$. To make sure that $I'$ satisfies the same conditions in Setup~\ref{setup}, we say that $I'$ is \emph{similar} to $I$ if $S(I')=S(I)$, $I'$ also dominates $R \cup S(I)$, and $F_{I'}$ has the same components as $F_{I}$. It follows that $|I'|=|I|=2(t-d-1)$, and, in particular, that $I'$ satisfies the assumptions of Lemma~\ref{prop b} and 
Lemma~\ref{avsize} \ref{old avsize}.

In the following lemma, Parts (ii) and (iii) resemble \cite[Lemma 3.1 (i) and (ii)]{HS06}, respectively. Their proofs are similar to the ones in \cite{HS06} but require more work because our IMC has $d+1$ components and only dominates $t\le r$ classes (in contrast, IMCs in \cite{HS06} have only one component and dominate the vertices in all $r$ classes of $G$). 
Furthermore, in Lemma~\ref{samecomp} (iii), we need to show that the new IMC is similar to $I$.
%In particular, we need to apply Lemma~\ref{prop b} to extend a partial IT on $R\cup S(I)$ to a partial IT of $G$.
\begin{lemma}\label{samecomp}
Let $G$ be as in Setup~\ref{setup}~(i). 
%all the vertices in $R\cup S(I)$, where $|R\cup S(I)|=t$.
Suppose $v, w \in I$ are adjacent and contained in the component $C$ of $F_I$. Then the following holds.
\begin{enumerate}[label=\textup{(\roman*)}]
\item All vertices in $A_v\cup A_w$ are in $C$. %the same component as $v$ and $w$.
%\item Every vertex in $A_v$ is in the same component as $w$ upon deleting the edge $\{v,w\}$.
\item $G[A_v,A_w]$ is a complete bipartite graph. 
\item For any $a\in A_v$ and $b\in A_w$, $(I\setminus\{v, w\}) \cup \{a, b\}$ is an IMC in $G$ similar to $I$.
%all the vertices in $R\cup S(I)$.
\end{enumerate} 
\end{lemma}
\begin{proof}
For a vertex $x\in V(G)$, let $V(x)$ denote the class $V_i\ni x$.

 (i) Suppose some vertex $a\in A_v$ is contained in $C'$ for some component $C'\ne C$ of $F_I$. By Lemma \ref{it from imc}, in every component of $F_I$, we can find a subset of $I$ forming a partial IT covering all but one arbitrary class.  
 We pick a partial IT in $C$ omitting the class of $v$, a partial IT in $C'$ omitting the class of $a$, and partial IT's omitting an arbitrary class in the other components.  Let $T_1$ be the union of these partial IT's. By Lemma~\ref{prop b}, there is an IT $T_2$ of the classes not dominated by $I$ such that $T_2\cap N_G(I\cup \{a\}= \emptyset$. Then $T_1\cup T_2 \cup \{a\}$ is an $(r-d)$-IT of $G$ because $N(a)\cap (I\setminus \{v\}) = \emptyset$.

%Observe that $a \not\in N(w)$ as otherwise $a$ would not be uniquely dominated by $v$. there exists a partial IT $T$ containing a vertex in $I$ from all but one of the classes in each component of $I$. We may assume that $T$ contains no vertex from $V(v)$ and no vertex from $V(a)$. Then, by Lemma~\ref{prop b}, $T$ can extend to $U$ to form a partial IT $T'$ of size $(r-d-1)$ in $G$, independent on $a$. But then $T' \cup {a}$ would be an IT of size $(r-d)$ in $G$, a contradiction. 

%A proof of Part 2 can be found in the proof of \cite[Lemma 3.1 (i)]{HS06}. 
\medskip
(ii) Consider $a \in A_v$ and $b \in A_w$. 
Part (i) shows that $a, b$ are both in $C$.
Let $F_I - vw$ be the forest obtained by removing the edge $vw$ from $F_I$, and let $C_v, C_w$ be two components of $F_I - vw$ that contain $v, w$, respectively. We claim that
\begin{align}\label{eq:aC}
  V(a) \in C_w, \quad \text{and} \quad V(b) \in C_v.  
\end{align}
%We need the following claim.
%\begin{claim}
%For every $a \in A_v, b \in A_w$, we have $a \in C_w, b \in C_v$.
%\end{claim}
%\begin{proof}[Proof of Claim]
%We only show for $a \in C_w$ as the proof is similar for $b$. 
%Accordingly, $I - \{v, w\}$ becomes two IMC's $I_v$ and $I_w$ on $C_v, C_w$, respectively. 
Indeed, suppose that $V(a) \not\in C_w$. Then $V(a)\in C_v$. Let $I_v= (I\setminus \{v\}) \cap \bigcup_{V_i\in C_v} V_i$. Then $I_v$ is an IMC of $G[\bigcup_{V_i\in C_v} V_i]$. By Lemma \ref{it from imc}, we can find a set $T_1\subset I\cap C_v$ as a partial IT avoiding the class of $a$. Similarly we find a partial IT $T_2$ in $C_w$ avoiding the class of $w$. The union $T_1\cup T_2\cup \{a, w\}$ is a full $IT$ of $C$ because $N(a)\cap (I\setminus \{v\}) = \emptyset$.  We then extend it to an $(r-d)$-IT of $G$ by Lemma~\ref{prop b}, which contradicts our assumption. The same arguments show that $V(b) \in C_v$.

%\end{proof}
%For (ii), we first show that $a \in C'_w$. If $a \in C'_v$ then since $a$ is not dominated by $I-\{v\}$, $I$ induces an IT $T_a$ in $C'_v$ containing $a$. And as $w$ is not dominated by the IMC induced in $C'_w$, $I$ induces an IT $T_w$ in $C'_w$ containing $w$. Now, since $a$ is independent of $w$, $T_a \cup T_w$ gives an IT $T$ in $C$. By Lemma~\ref{prop b}, we can extend $T$ to $U$ which, combined with the induced partial IT of maximum size in all other components of $I$, gives an $(r-d)$-IT in $G$, a contradiction. Thus $a \in C'_w$. 

%We claim that ${a,b} \in E(G)$. Consider the component $C$ of $I$ containing $V(v)$ and $V(w)$. By the first part, both $A_v$ and $A_w$ are contained in this component. Consider the graph $C'$ obtained by removing $\{v,w\}$, i.e. setting $V(C')=V(C)-\{v,w\}$ and $E(C')=E(C)-\{v,w\}$. Let $C'_v$ be the component of $C'$ containing the class of $v$ and $C'_w$ be the component containing the class of $w$. Observe that $I$ induces an IMC $I_v$ in $C'_v$ and an IMC $I_w$ in $C'_w$. By the first part, $a \in C'_w$ and $b \in C'_v$. 

Now suppose that $ab \not\in E(G)$. Then, we can obtain partial ITs $T_1\subseteq I_w$ and $T_2\subseteq I_v$ 
such that $T_1$ covers all classes of $C_w$ except for the class of $a$ and $T_2$ covers all classes of $C_v$ except for the class of $b$. Since $ab \not\in E(G)$, $N(a)\cap (I\setminus \{v\}) = \emptyset$, and $N(b)\cap (I\setminus \{w\}) = \emptyset$, the union $T_1 \cup T_2 \cup \{a,b\}$ is a full IT of $C$, which gives an IT of size $r-d$ in $G$ after being extended to $U$ by Lemma~\ref{prop b} , a contradiction. This shows that $G[A_v,A_w]$ is a complete bipartite graph. 

\medskip
(iii) We first show that $I':= (I\cup\{a\})\setminus \{w\}$ is an IMC of $G$ similar to $I$.
Indeed, since $I$ is an IMC and $N(a)\cap I' = \{v\}$, 
%$a$ only connects to $v$ and $v$ only connects to $a$ in $I'$, 
$I'$ is an induced matching in $G$. Recall that $V(a) \in C_w$ and $V(v) \in C_v$ by \eqref{eq:aC}, where $C_v, C_w$ are the two components of $F_I - vw$ that contain $v, w$, respectively.
Since the edge $av$ connects $C_v$ and $C_w$, the graph $\G_{I'}$ is a forest with the same components as $\G_{I}$. 
%Thus the induced graph of $I'$ is a connection of $C_w$ and $C_v$ and hence a forest on $S(I) \cup R$ with $d+1$ components. 
We claim that $I'$ dominates all the classes of $R\cup S(I)$. Indeed, let $x \in \bigcup_{V_i\in R\cup S(I)} V_i$. If $x$ has a neighbor in $I\setminus \{w\}$, then $x$ is dominated by $I'$. Otherwise $x$ is only adjacent to $w$ in $I$, giving $x \in A_w$ and thus $ax \in E(G)$ by Part (ii). Hence, $x$ is dominated by $I'$ as desired. 

We now show that $S(I')=S(I)$.
To see $S(I') \subseteq S(I)$, it suffices to show that $V(a) \in S(I)$. By (i), $a$ is in $C$ so $V(a) \cap I \neq \emptyset$. To see $S(I) \subseteq S(I')$, it suffices to have $V(w) \in S(I')$. If this is not the case, then $V(w)\ne V(a)$ and $V(w)$ is a leaf in $F_I$. 
Then $V(a)\notin C_w = \{V(w)\}$, contradicting \eqref{eq:aC}. 

%Thus $I'$ is an IMC in $G$ with the same components as $I$ and dominating all the vertices in $R\cup S(I)$. since $a\in A_v$, it follows that

For $u\in I'$, recall $A_u(I') = \{x\in \bigcup_{V_i\in R\cup S(I)} V_i: N(x)\cap I' = \{u\}\}$. 
We have $b\in A_a(I')$ because $N(b)\cap I = \{w\}$ and $ab\in E(G)$. The arguments in the previous paragraph show that $I' \cup \{b\} \setminus \{v\} = (I\setminus\{v, w\}) \cup \{a, b\}$ is an IMC in $G$ similar to $I$. This completes the proof of (iii).
%Now, note that $a \in A'_v$ and $b \in A'_a$. Thus by the previous paragraph $(I\setminus\{v, w\}) \cup \{a, b\}=I' \cup \{b\} \setminus \{v\}$ is an IMC that dominates all of the vertices in $R\cup S(I)$. 
\end{proof}
Suppose that $I$ is an IMC on $2(t-d+1)$ vertices. 
Let $H$ be the $2(t-d-1)$-partite graph on $\{A_v \colon v \in I\}$ obtained from $G$ by removing edges in $G[A_v]$ for all $v\in I$ and edges between $A_v$ and $A_w$ for all $vw\in E(G[I])$. \emph{An $IT$ of} $H$ is an independent set of $H$ with one vertex from each $A_v$ (not from original classes $V_i$ of $G$), for example, $I$ is an IT of $H$. 

The following lemma and proof are similar to \cite[Lemma 3.3]{HS06} and its proof. 
\begin{lemma}\label{it is imc}
Let $G$ be as in Setup~\ref{setup}~(i). Every IT $T$ of $H$ is an IMC in $G$ similar to $I$.    
\end{lemma}
\begin{proof}
Suppose that $I=\{v_i, w_i: 1\le i\le t-d-1\}$ with $v_i w_i\in E(G)$ and $T= \{a_i, b_i: 1\le i\le t-d-1\}$ with $a_i \in A_{v_i}$ and $b_i \in A_{w_i}$. 
For $1\le j \le t-d-1$, we claim that 
\[
I_j := \left(I\setminus \bigcup_{i\le j} \{v_i, w_i\} \right) \cup \bigcup_{i\le j} \{a_i, b_i\}
\]
is an IMC in $G$ similar to $I$. In particular, $T=I_{t-d-1}$ is similar to $I$, as desired.

We prove the claim by the induction on $j$. The $j=1$ case is 
Lemma \ref{samecomp} (iii). Suppose $j\ge 1$ and $I_j$ is similar to $I$. We note that $I_{j+1}= (I_j \setminus \{v_{j+1}, w_{j+1}\}) \cup \{a_{j+1}, b_{j+1}\}$. For simplicity, we write $v, w, a, b$ by omitting subscripts ${j+1}$. We know $a\in A_v(I)$ and thus $N(a)\cap I= \{v\}$. Since $T$ is an independent set, it follows that $N(a)\cap I_j = \{v\}$. Thus $a\in A_v(I_j)$. Similarly $b\in A_w(I_j)$. We can thus apply Lemma \ref{samecomp} (iii) and conclude that $I_{j+1}$ is an IMC in $G$ similar to $I$. 
%$I' := I \cup \{a_1, b_1\} \setminus \{v_1, w_1\}$ is an IMC that dominates all of the vertices in $R \cup S(I)$. Then, for all $i>1$ we still have $a_i \in A'_{v_i}$ and $b_i \in A'_{w_i}$ since $T$ is an IT of $H$. Therefore, by repeating this argument we can see that $T$ is an IMC that dominates all of the vertices in $R \cup S(I)$.
\end{proof}

\begin{lemma}\label{lem a}
Suppose that $G$ is as in Setup~\ref{setup}~(ii) and let $Y \subseteq I.$ If $A'_v \subseteq A_v$ for $v \in Y$ are subsets such that $A'_v \neq \emptyset$ for all $v \in Y$ and $\sum_{v \in Y}|A_v \setminus A'_v| \leq\Delta$, there exists an IT of $\{A'_v \colon v \in Y\}$ in $H$. 
%Thus, for any vertex $x$ in $G$ such that $A_v \setminus N(x) \not=\emptyset$ for all $v \in I,$ we can see by setting $A'_v=A_v\setminus N(x)$, there exists an IT of $H$ independent of $x$.   
\end{lemma}
\begin{proof}
Suppose that $\{A'_v \colon v \in Y\}$ contains no IT in $H$. Then, by Lemma~\ref{lem:noit}, there exists a set $S \subseteq Y$ and an nonempty set $Z$ of edges on $V_S :=\{A'_v \colon v \in S\}$ such that (1) $X:=V(Z)$ dominates $V_S$ in $H$, (2) $|X| \le 2(|S|-1)$, and (3) $X\cap A'_v\ne \emptyset$ for every $v\in S$.
%\begin{itemize}
%\item $X:=V(Z)$ dominates $V_S$ in $H$;
%\item $|X| \le 2(|S|-1)$;
%\item $X\cap A'_v\ne \emptyset$ for every $v\in S$.
%\end{itemize}

Since $X$ dominates $V_S$ in $H$ and $\sum_{v \in Y}|A_v \setminus A'_v| \leq \Delta$, we have 
\begin{align}\label{eq:39}
\Delta |X| - \sum_{x \in X}\deg_{G\setminus H}(x) \ge \sum_{x \in X}\deg_H(x) \ge \sum_{v \in S}|A'_v| \ge \sum_{v \in S}|A_v|-\Delta.
\end{align}
For every $x\in X$, by the definition of $H$ and Lemma \ref{samecomp}, we have $\deg_{G\setminus H}(x) \ge |A_{w_x}|$, where $w_x$ is the neighbor of $v\in I$ in $G[I]$ such that $x\in A_v$. Since $X\cap A'_v\ne \emptyset$ for every $v\in S$, there exists $X'\subseteq X$ such that $|X'\cap A'_v|=1$ for every $v\in S$. Since each $v\in S$ is adjacent to a unique $w\in I$, we have   
\[
 \sum_{x \in X}\deg_{G\setminus H}(x) \ge \sum_{x \in X} |A_{w_x}| 
 \ge \sum_{x \in X'} |A_{w_x}|  > (|S| -1 )\Delta 
 \]
by Lemma \ref{avsize} \ref{old avsize}.  Together with \eqref{eq:39} and the assumption $|X|\le 2(|S| -1)$ , we obtain that
\begin{align}\label{eq:DS}
 \Delta(2|S|-1)\ge \Delta(|X|+1)|\ \ge \sum_{x \in X}|A_{w_x}|+\sum_{v \in S}|A_v|
> 2(|S| -1) \Delta   
\end{align}
by applying Lemma \ref{avsize} \ref{old avsize} again.
%We now show that $\sum_{x \in X}|A_{w_x}|+\sum_{v \in S}|A_v|>(2|S|-1)\Delta$ to give the desired contradiction.
%Observe that guarantees that $\sum_{x \in X}|A_{w_x}|+\sum_{v \in S}|A_v >(2|S|-2)\Delta$. 

If $|X|<2(|S|-1)$ then $\Delta(|X|+1)\leq(2|S|-2)\Delta$, giving a contradiction. We thus assume that $|X|=2(|S|-1)$. %We also assume that $X$ is the vertex set of $(|S|-1)$ edges, the set of which we denote by $Z$. 
Let $T$ be the tree obtained from $(V_S, Z)$ by contracting each $A'_v$ into a vertex. 
%$(S,E)$ where we define $E$ such that $\{v_1, v_2\} \in E$ if $\{x_1, x_2\} \in Z$ for some $x_1\in A_{v_1}, x_2 \in A_{v_2}$. 
We proceed by cases on $q := |\{v \in S \colon |A_v \cap X|\ge 2\}$.

Suppose $q=0$. Then, $|Z|=1$ and $|S|=|X|=2$. Suppose $S=\{v_0, v_1\}$. We must have ${v_0, v_1} \not\in E(G[I])$ as otherwise there would be no edges connecting $A_{v_0}$ and $A_{v_1}$ in $H$. Thus $v_0,w_0,v_1,w_1$ are all distinct and we may calculate
\[\sum_{v \in S}|A_v|+\sum_{x \in X}|A_{w_x}|=|A_{v_0}|+|A_{w_0}|+|A_{v_1}|+|A_{w_1}|>3\Delta,\]
where the last inequality follows from Lemma \ref{avsize} \ref{old avsize}. 

Now consider the case $q=1$. Then $T$ is a star with at least two edges. Moreover, let $s = |S|\ge 3$. Let $v_1$ be the distinct vertex in $S$ with $|A_v \cap X| \ge 2$. Suppose $x_1,x_2,\ldots,x_{s_1} \in X$. Then $Z=\{\{x_1,y_1\},\{x_2,y_2\},\ldots,\{x_{s_1},y_{s-1}\}\}$, where $y_1,\ldots,y_{s_1}$ are in distinct members of $\{A_{v_2},A_{v_3},\ldots,A_{v_s}\}$. As before, $v_1$ is not adjacent to any other $v_i$ in $I$ as this would imply that there are no edges between $A_{v_1}$ and $A_{v_i}$ in $H$. So,
\begin{align*}
\sum_{v \in S}|A_v|+\sum_{x \in X}|A_{w_x}|&=\sum_{v \in S}|A_v|+\left(\sum_{w \in V(I) \colon N(w)\cap S \neq \emptyset}|A_w|\right)+|A_{w_1}|(s-2)\\
&>\sum_{v \in S\cup \{w_1\}}|A_v|+\sum_{w \in V(I) \colon N(w)\cap S \neq \emptyset} |A_w|\\
&>s\Delta+(s-1)\Delta,
\end{align*}
where again the last inequality follows from Lemma \ref{avsize} \ref{old avsize}.

Finally, let $q \ge 2$. Let $v_1,v_2$ be distinct vertices in $S$ such that $|A_{v_i}\cap X|\ge 2$ for $i=1,2$. Then,
\begin{align*}
\sum_{v \in S}|A_v|+\sum_{x \in X}|A_{w_x}|&\ge \sum_{v \in S}|A_v|+\left(\sum_{w \in V(I) \colon N(w)\cap S \neq \emptyset}|A_w|\right)+|A_{w_1}|+|A_{w_2}|\\
&>(s-1)\Delta+(s-1)\Delta+\Delta.\qedhere
\end{align*}
In each case we derived that $\sum_{v \in S}|A_v|+\sum_{x \in X}|A_{w_x}|> (2|S|-1)\Delta$, which contradicts \eqref{eq:DS}.
\end{proof}

We now derive a key property on the structure of $G$.

\begin{lemma}\label{notdomsamecomp}
Let $G$ be as in Setup~\ref{setup}~(ii). Then, every vertex $x\in \bigcup_{V_i\in R\cup S(I)} V_i$ is completely joined to some $A_v$ where $v\in I$ is in the same component as $x$.
%and adjacent to some vertex of $I$ in the same  component of $F_I$. 
\end{lemma}
\begin{proof}
Fix $x\in \bigcup_{V_i\in R\cup S(I)} V_i$.
If $x\in A_w$ for some $w\in I$ such that $vw\in E(G[I])$, 
then by Lemma \ref{samecomp}, $x$ is completely joined to $A_v$ and all the vertices in $A_v$ are the same component as $x$. 

We thus assume that $x\not\in \bigcup_{v\in I} A_v$. We first show that $x$ is completely joined to some $A_v$.
Suppose, to the contrary, that $A_v\setminus N(x) \ne \emptyset$ for all $v \in I$. Applying Lemma~\ref{lem a} with $Y=I$ and $A'_v = A_v\setminus N(x)$, we obtain an IT $I'\subset \bigcup _{v\in I} A'_v$ of $H$ that is independent of $x$. Now, Lemma \ref{it is imc} implies that $I'$ is an IMC in $G$ that dominates all the vertices of $R\cup S(I)$, including $x$. This contradicts the fact that $I'\cap N(x) = \emptyset$.
%By Lemma~\ref{prop b}, there is a partial IT $T'$ of $U$ that is independent of $T\cup \{x\}$. By Lemma \ref{it from imc}, there is an $(r-d-1)$-IT $T''\subset T$ avoiding the class of $x$ . The union $T'\cup T''\cup \{x\}$ is an $(r-d)$-IT of $G$, giving a contradiction. 

Suppose $x$ is completely joined to $A_{v_1}$ for some $v_1\in I$ but $x$ and $A_{v_1}$ are in different components.  
If $x$ is completely joined to $A_{v'}$ for some $v'\in I\setminus \{v_1\}$, then, by Lemma \ref{avsize} \ref{old avsize}, $\deg_G(x)\ge |A_{v_1}|+|A_{v'}|> \Delta$, a contradiction. Thus, $A_{v}\setminus N(x)\ne \emptyset$ for all $v\in I\setminus \{v_1\}$.
Applying Lemma \ref{lem a} with $A'_{v_1}=A_{v_1}$ and $A'_v=A_v-N(x)$ for $v \in I \setminus {v_1}$ (thus $A'_v\ne \emptyset$), we obtain an IT $I'\subset \bigcup _{v\in I} A'_v$ of $H$. 
By Lemma \ref{it is imc}, $I'$ is an IMC in $G$ that dominates all the vertices in $R \cup S(I)$ and has the same components as $I$.
Suppose $\{u\} = I'\cap A'_{v_1}$. Then $x\in A_u(I')$ because $N(x)\cap I' = \{u\}$. However, $x$ and $u$ are in different components, contradicting Lemma \ref{samecomp} (i).
\end{proof}

We are now ready to prove Theorem~\ref{genbound}.
\begin{proof}[Proof of Theorem~\ref{genbound}]
Suppose $G$ has no $(r-d)$-IT. 
Add edges if necessary so that the largest partial IT of $G$ has size $r-d-1$.
Let $I$ be the IMC given by Theorem \ref{imc exists}. Let $J$ be the smallest component of $F_I$ with $j\le q$ classes. By Lemma \ref{notdomsamecomp}, all the vertices in $V_J:=\bigcup_{V_i\in J} V_i$ are dominated by some vertex in $I\cap V_J$. Since $F_I[J]$ is a tree on $j$ vertices and $I$ is an IMC, it follows that $|I\cap V_J| = 2(j-1)$. However, since 
\[
j n-2\Delta(j-1)>2j\Delta\left(1-\frac{1}{q}\right)-2\Delta(j-1)=2\Delta\left(1-\frac{j}{q}\right)\ge 0,
\]
$2(j-1)$ vertices cannot dominate $j n$ vertices in $G$, giving the desired contradiction.
\end{proof}

\section{Proofs of Theorems~\ref{oddq} and~\ref{odd q}}
\label{sec:odd q results}
Let $G$ be an $r$-partite graph, where $r=q(d+1)+k$, $0\le k\le d$, and $q$ is odd.
We have a new setup due to the new lower bounds for $n$ in Theorems~\ref{oddq} and~\ref{odd q}.
%As before, let $r=q(d+1)+k$, where $k<d+1$.  

\begin{setup}\label{odd setup}
We have the following two assumptions:
\begin{enumerate}[label=\textup{(\roman*)}]
    \item Let $G$ be as in Setup~\ref{setup}~(ii) with the additional assumption that $n>2\Delta\left(1-\frac{1}{q-1}\right)$.
    \item Let $G$ be as in Setup~\ref{odd setup}~(i) with the additional assumption $n>\Delta\left(2-\frac{6d+7}{3r}\right)$.
    %$n>2\Delta-(4d+5)\Delta/2r$.
\end{enumerate}
\end{setup}

Since $2-\frac{6d+7}{3r}>2-\frac{4d+5}{2r}$, Setup~\ref{odd setup} (ii) has the largest bound for $n$. Furthermore, when $G$ is as in Setup~\ref{odd setup}, we can use all the lemmas in Section~\ref{sec:IMC}. 
%For convenience, let $V'= \bigcup \{V_i: V_i\in S(I)\cup R\}$.  Then $|V'|= tn$ and $|I|= 2(t-d-1)$ for some $t\le r$.

We first show that if $G$ is as in Setup~\ref{odd setup}~(i), then every component of $F_I$ has at least $q$ vertices.
%that is, $R\cup S(I)=\{V_1, \dots, V_r\}$.  
\begin{lemma}\label{t=r}
Let $G$ be as in Setup~\ref{odd setup}~(i). Then the forest $F_I$ has at least $q(d+1)$ vertices and exactly $d+1$ components, each of which has at least $q$ vertices.
\end{lemma}
%that Then for every part $V_i$ of $G, V_i \in S(I)$.
%with vertex partition $V_1 \cup \cdots \cup V_r$ 
\begin{proof}
Since $F_I$ has $d+1$ components, 
it suffices to show that every component $J$ of $F_I$ has at least $q$ vertices. Suppose that $J$ has $b$ vertices. 
Then $|I\cap J| = 2(b-1)$ because $\G_{I}[J]$ is a tree.
By Lemma \ref{notdomsamecomp}, every vertex of $G$ in $J$
is dominated by a vertex in $I\cap J$. Thus,
\[2\Delta(b-1)\ge bn > b2\Delta\left(1-\frac{1}{q-1}\right),\]
which implies that $b>q-1$. Since $b$ is an integer, it follows that $b \geq q$.
\end{proof} 

We can now prove Theorem~\ref{oddq}. The main idea is to find a component of $F_I$ with $q$ vertices and apply \eqref{eq:nr} to find a full IT in this component by using the vertices of $\bigcup_{v\in I} A_v$.

\begin{proof}[Proof of Theorem~\ref{oddq}]
Suppose to the contrary, that $G$ has no $(r-d)$-IT. Since
$n>\Delta(2-\frac{4d+5}{2r})$, by Theorem~\ref{imc exists}, $G$ has an IMC $I$ of $d+1$ components that dominates $S(I)\cup R$. Since $n> 2\Delta(1- \frac{q}{q^2-1}) > 2\Delta (1 - \frac{1}{q-1})$, Lemma \ref{t=r} says that $F_I$ has at least $q(d+1)$ vertices and every component of $F_I$ has at least $q$ classes. Since $r=q(d+1)$, $F_I$ has exactly $r=q(d+1)$ vertices, and every component of $F_I$ has $q$ classes.

Let $D$ denote the set of vertices in $V(G)$ that are not in any $A_v$. By Lemma~\ref{avsize}~\ref{odd q max D}, $|D|\le 2\Delta(r-d-1)-rn$. 
We know that some component $J$ of $I$ contains at most $|D|/(d+1)$ vertices of $D$. Let $G'_J = G[V_J\setminus D]$ be the induced subgraph of $G$ formed by removing the vertices of $D$ from the vertex classes of $J$. Then each class of $G_J$ has size at least $n'$ with
\begin{align*}
n' &\geq n-\frac{|D|}{d+1} \geq n-\frac{2\Delta(r-d-1)-rn}{d+1}\\
&=n-\frac{2\Delta(q-1)(d+1)-q(d+1)n}{d+1}=n(q+1)-2\Delta(q-1)\\
&> 2\Delta(q+1)\left(1-\frac{q}{q^2-1}\right)-2\Delta(q-1)=2\Delta\left(1-\frac{1}{q-1}\right).
\end{align*}

By \eqref{eq:nr}, since $q$ is odd, the graph $G_J$
has a full IT $T_J$ consisting only of vertices of $V_J\setminus D$. We know $V_J\setminus D = \bigcup_{v\in J} A_v$ from Lemma~\ref{samecomp}. Hence, $T_J$ is independent of $I\setminus J$. Since for each component
$J'\not=J$, we can find in $J'$ an IT of all except one class of $J'$, these together with $T_J$ form an $(r-d)$-IT of $G$. This contradiction shows that $G$ indeed has an $(r-d)$-IT, and completes the proof.
\end{proof}

%When $r=6$ and $d=1$, Theorem~\ref{oddq} implies that $n>5\Delta/4$ guarantees an $5$-IT in any 6-partite graph $G$ with $n$ vertices in each part. This is tight due to Corollary \ref{main construction} with $q=2$ and $i=2$. This example also shows that Proposition~\ref{big construction}~\ref{construction iii} is tight. 

\bigskip
In order to prove Theorem~\ref{odd q}, we follow the approach of \cite{HS06} and obtain the structure of $G$.
To this end, we prove analogues of \cite[Lemma 3.5, Theorem 3.6, Theorem 3.7]{HS06} under Setup~\ref{odd setup}~(ii), in which we assume $n>\Delta\left(2-\frac{6d+7}{3r}\right)$.
%The statements and proofs of the generalizations of Lemma 3.5, Theorem 3.6, and Theorem 3.7 are as follows. 
\begin{lemma}\cite[Lemma 3.5]{HS06}
\label{almost bpg}
Let $G$ be as in Setup~\ref{odd setup}~(ii) and $V'= \bigcup \{V_i: V_i\in S(I)\cup R\}$. Suppose $a,b \in V'$ are completely connected to $A_w$ and $A_v$ respectively, where $vw \in E(G[I])$. Then $ab\in E(G)$. 
\end{lemma}
\begin{proof}
Suppose that $a$ is not adjacent to $b$. Then, for $u \in I-\{v,w\}$, let $A'_u=A_u \setminus N(a,b)$. If there exists $u \in I-\{v,w\}$ such that $A'_u$ is empty, then $\{a,b\}$ dominates $A_v,A_w$, and $A_u$, which contains more than $2\Delta$ vertices by Lemma \ref{avsize} \ref{old avsize}, contradicting $\Delta(G)\le \Delta$.  

Now, since $\{a,b\}$ dominates $A_v \cup A_w$, which has size more than $\Delta$ by Lemma \ref{avsize} \ref{old avsize}, we can see that $\sum_{u \in (I-\{v,w\})}|A_u\setminus A'_u| \leq\Delta$ and thus an IT $I'_0$ of $\{A'_u \colon u \in I-\{v,w\}\}$ exists by Lemma \ref{lem a}. Then $I':= I'_0 \cup \{v,w\}$ is an IT of $\{A'_u \colon u \in I\}$, which is an IMC of $G$ similar to $I$ by Lemma \ref{it is imc}.

By definition, we have $N(a,b)\cap I'=\{v,w\}$. Also, neither $a$ nor $b$ can individually dominate both $A_v(I')$ and $A_w(I')$ as they have size more than $\Delta$ combined by Lemma \ref{avsize} \ref{old avsize}. We now claim that there exist $v' \in A_v(I')$ and $w' \in A_w(I')$ such that $v'$ is not adjacent to $a$ and $w'$ is not adjacent to $b$. Then, $I'':= I'-\{v,w\}\cup\{v',w'\}$ is an IMC similar to $I$ by Lemma \ref{it is imc} such that $a \in A_{w'}(I'')$ and $b \in A_{v'}(I'')$ and therefore $ab \in E(G)$ by Lemma~\ref{samecomp}~(ii). 

Suppose that this is not possible. Assume $a$ dominates $A_{v}(I')$. By definitions, $A_{v}(I')\cap A_w(I) = \emptyset$, and the vertices in $A_v(I) \setminus A_v(I')$ are dominated twice in $I'$. Hence, we have 
\begin{align*}
\Delta &\ge \deg(a)\ge |A_w|+|A_v(I')| \ge |A_w|+|A_v|-|A_v \setminus A_v(I')|\\
&\ge 2tn - \Delta(4t - 4d - 6) - |A_v \setminus A_v(I')| &\text{by Lemma \ref{avsize}~\ref{odd q min Av}}\\
&\ge 3tn-\Delta(6t-6d-8) &\text{by Lemma \ref{avsize}~\ref{odd q max D}}.
\end{align*}
This implies $n\leq 2\Delta-\frac{6d+7}{3t}\Delta \leq2\Delta-\frac{6d+7}{3r}\Delta$, a contradiction. 
\end{proof}
\begin{lemma}\cite[Lemma 3.6]{HS06}
\label{Lem:3.6}
Let $G$ be a graph with a vertex partition $V_1\cup \cdots \cup V_r$ 
such that $|V_i|> 2\Delta(1- \frac{2d+3}{2r})$. Suppose $G$ has no 
$(r-d)$-IT but $G-e$ has an $(r-d)$-IT for some edge $e\in E(G)$. Then $e$ lies in an IMC of $G$ returned by Theorem \ref{imc exists}.
\end{lemma}
\begin{proof}
Let $e=\{v_1,v_2\}$. Then, by assumption there exists a transversal $T'=\{v_1,v_2, \ldots, $ $v_{r-d}\}$ such that $T'-\{v_j\}$ is independent for $j=1,2$. Let $I_0=\{v_1,v_2\}$ and $T_0=T'-\{v_1\}$. It is easy to check that $(I_0, T_0)$ is a feasible pair in $G$. Indeed, Conditions~(1)--(4) of Definition~\ref{feasible pair} are trivially satisfied, and Condition~(5) follows from the observation that $W=\{v_1\}$ so any $T$ contradicting Condition~(5) is an $(r-d)$-IT in $G$. Then, Theorem~\ref{imc exists} with the same $(I_0, T_0)$ gives an IMC $I$ in $G$ which contains $e$ by Theorem~\ref{imc exists}~(i).
\end{proof}

\begin{lemma}\cite[Lemma 3.7]{HS06}
\label{full bpg}
Let $G$ be as in Setup \ref{odd setup}~(ii). Moreover, assume that $G-e$ has an $(r-d)$-IT for every $e \in G$. Let $V'= \bigcup \{V_i: V_i\in S(I) \cup R\}$. Then $G[V']$ is a union of $t-d-1$ vertex-dsijoint complete bipartite graphs. 
\end{lemma}

\begin{proof}
%As in the proof of Lemma~\ref{it is imc}, 
We label $I=\{v_i, w_i: 1\le i\le t-d-1\}$ with $v_i w_i\in E(G)$.
By Lemma~\ref{notdomsamecomp}, every vertex $x\in V'$ is completely joined to some $A_v$.  
By Lemma~\ref{avsize} \ref{old avsize}, such $A_v$ is unique (otherwise $\deg(x)> \Delta)$.
By Lemma~\ref{almost bpg}, we can write $G[V']$ as a disjoint union $A_1 \cup \cdots A_{t-d-1} \cup B_1 \cup \cdots \cup B_{t-d-1}$ such that $G[A_i,B_i]$ is complete, $A_{v_i}\subseteq A_i$, and $A_{w_i} \subseteq B_i$ for all $i$. Thus, it remains to show that there are no edges outside these bipartite subgraphs.

Suppose $e=xy\notin \bigcup_{1 \leq i \leq t-d-1} E(G[A_i, B_i])$. By Lemma \ref{Lem:3.6}, there is an IMC $I'$ in $G$ returned by Theorem~\ref{imc exists} containing $e$. Without loss of generality, let us write $x \in A_1$. Then by assumption, $y \not\in B_1$. Suppose that $y \in A_1$. Then $B_1$ is dominated by $\{x, y\}\subseteq I'$. Thus $|B_1|\leq2\Delta(t-d-1)-tn$ by applying Lemma \ref{avsize}~\ref{odd q max D} with $I'$. On the other hand, by applying Lemma \ref{avsize}~\ref{odd q min Av} with $I$, we get $|B_1|\ge 2tn-\Delta(4t-4d-5)$. Together, these two inequalities give $n \leq 2\Delta-\frac{6d+7}{3t}\Delta \leq2\Delta-\Delta\frac{6d+7}{3r}$, a contradiction. 

Thus we may assume $y \not\in A_1 \cup B_1$. Without loss of generality, assume $y \in A_2$. Considering the IMC $I'$, we have that $x \in A_y(I')$ and $y \in A_x(I')$. Suppose that $A_y(I') \cap B_2$ or $A_x(I') \cap B_1$ is empty, for example, $A_y(I') \cap B_2=\emptyset$. 
Since $B_2$ is dominated by $y$, it follows that $B_2$ is dominated at least twice in $I'$ and the argument in the previous paragraph gives the same contradiction.

We are left with the case that both $A_y(I') \cap B_2$ and $A_x(I') \cap B_1$ are nonempty. Let $u \in A_y(I') \cap B_2$ and $w \in A_x(I') \cap B_1$. By Lemma~\ref{samecomp}, we know $A_x(I')\subseteq N(u)$ and $A_y(I')\subseteq N(w)$. Consequently,  
\begin{align*}
2\Delta &\ge \deg(w) + \deg(u) \ge (|A_1|+|A_y(I') \cap B_2|)+(|A_2|+|A_x(I') \cap B_1|)\\
&\ge (|A_1|+|A_y(I')|+|B_2|-|A_y(I') \cup B_2|)+(|A_2|+|A_x(I')|+|B_1|-|A_x(I') \cup B_1|)\\
&\ge (|A_1|+|A_2|+|B_1|+|B_2|)+(|A_x(I')|+|A_y(I')|)-(|A_y(I') \cup B_2|+|A_x(I') \cup B_1|)\\
&> 3\Delta+\Delta-(|A_y(I') \cup B_2|+|A_x(I') \cup B_1|) \qquad \text{(by Lemma \ref{avsize} \ref{old avsize})}.
\end{align*}

Since $A_x(I') \cup B_1$ and $A_y(I') \cup B_2$ are dominated by $x$ and $y$ respectively, each has size at most $\Delta$. Thus, $2\Delta>4\Delta-\Delta-\Delta=2\Delta$, a contradiction. 
\end{proof}

We are ready to prove Theorem~\ref{odd q}. The main idea is, due to Lemmas~\ref{samecomp} and~\ref{full bpg}, the components of $F_I$ give rise to a partition of $G[V']$ into $d+1$ components that are independent of each other. Once we apply \eqref{eq:nr} and find a full IT in one of the components, we immediately obtain an $(r-d)$-IT in $G$.

\begin{proof}[Proof of Theorem~\ref{odd q}]
%First, note that $n>2\Delta\left(1-\frac{d+1}{r-1}\right)\ge 2\Delta\left(1-\frac{1}{q-1}\right)$.

Suppose that $G$ has no $(r-d)$-IT. 
After removing edges from $G$ if necessary, we further assume that $G-e$ has an $(r-d)$-IT for every $e\in E(G)$.

Let $I$ be an IMC of $t-d+1$ edges given by Theorem \ref{imc exists}. Let $V'= \bigcup \{V_i: V_i\in S(I) \cup R\}$.
By Lemma~\ref{full bpg}, $G[V']$ is the union of $t-d-1$ vertex-disjoint complete bipartite graphs on $A_i \cup B_i$ for $i\leq t-d-1$. Since $G[I]$ is an induced matching of size $t-d+1$, each $A_i \cup B_i$ contains exactly one edge of $I$. Since there is no edge of $G$ between $A_i \cup B_i$ and $A_j \cup B_j$ for $i\ne j$, we have $A_v= A_i$ and $A_w = B_i$ if $v\in I\cap B_i$ and $w\in I\cap A_i$. 

%with $A_{v_i}\subseteq A_i$ and $A_{w_i}\subseteq B_i$ for all $i\leqr-d-1$. This means that every vertex is in some $A_v$, i.e. $D=\emptyset$. 
%This is because if some vertex $x$ is connected to two vertices $a,b \in I$, then there is an edge in $G$ outside of a complete bipartite graph, a contradiction. 
By Lemma~\ref{samecomp}~(i), all the vertices of $A_i\cup B_i$ lie in the same component of the forest $F_I$. Therefore, each component of $F_I$ is independent of other components because there are no edges between $A_i \cup B_i$ and $A_j \cup B_j$ for $i\ne j$. 

By Lemma~\ref{t=r}, every component of $F_I$ has at least $q$ classes of $G$. Since $q=\lfloor r/(d+1) \rfloor$, there exists a component $J$ of $F_I$ with exactly $q$ vertices. As $q$ is odd and $n>2\Delta\left(1-\frac{1}{q-1}\right)$, by \eqref{eq:nr}, there is a full IT in $J$. We can find by Lemma \ref{it from imc} an IT in every other component of $F_I$ missing a vertex from at one class in that component. Since each component is independent of the other components, combining these ITs gives an IT $T_0$ of size $t-d$ in $G$ contained in $S(I) \cup R$. By Lemma \ref{prop b}, there exists an IT in $G$ of size $(t-d)+(r-t)=r-d$ containing $T_0$, contradicting $G$ having no $(r-d)$-IT.
\end{proof}

\section{Concluding Remarks}

Let $r,d \geq 0$ and $q,k \geq 0$ be integers such that $r=q(d+1)+k$, where $k \leq d$. In this paper we have completely determined $n(r,d+1,\Delta)$ 
%when $q$ is even and $q \geq 4k$, and when $q$ is odd, $q \geq 6d+7$, and $k=0$ 
when $q \geq 4k$ is even and when $q \geq 6d+6k+7$ is odd. We have also shown that $n(6,2,\Delta)=\lfloor 5\Delta/4 \rfloor$, answering a specific question of \cite{LTZ22}. It is interesting to know the value of $n(r,d+1,\Delta)$ in the remaining cases:

\begin{enumerate}
\item $q$ is even and $q<4k$, such as $r=7$ and $d=2$,
\item $q$ is odd and $q<6d+6k+7$, such as $r=7$ or $10$ and $d=1$.
\end{enumerate}
The aforementioned results of \cite{LTZ22} determined $n(r,d+1,\Delta)$ in many cases when $d>r/2-1$ (equivalently, $q=1$) but there are still unknown ones such as $r=7$ and $d=3$.

One may also ask in an $r$-partite graph whether, given an $(r-d)$-IT $T$, there exists an $(r-d')$-IT $T'\supset T$ (provided $d'<d<r$). A similar proof to that of \cite[Corollary 15]{H16}, which uses topological methods, shows that if $G$ is an $r$-partite graph with $n$ vertices in each part, $0 \le d<r$, and $S$ is a $k$-IT in $G$ with $k<r-d$, then there exists an $(r-d)$-IT containing $S$ if 
\[
n>2\Delta\left(1-\frac{d+1-(k/2)}{r-k}\right).
\]
It would be interesting to know if any of the methods used in this paper, in particular, the results on IMCs developed in Sections \ref{sec:IMC} and \ref{sec:odd q results}, can be used to improve this bound.

Finally, it was shown in \cite{GH20} that for any $\Delta$, there exists an algorithm that, given a multipartite graph $G$ with at least $2\Delta+1$ vertices in each part, returns an independent transversal in polynomial time of $|V(G)|$. One may ask if there is a similar polynomial-time algorithm that returns $(r-d)$-ITs under a weaker condition on the size of parts.

%\section{Notes on odd $q$, $k>0$}

%Let $G$ be $r$-partite with no $(r-d)$-IT and $r=q(d+1)+k$ where $q$ is odd %and $k \geq 0$. Suppose that $n>2\Delta(1-1/(q-1))$. Let $I$ be an IMC %in $G$ with $d+1$ components.

%The same proof as Lemma 4.2 shows that each component of $I$ has at least $q$ %classes. This means that $|S(I)| \geq q(d+1)$, so $t \geq q(d+1)$. %It may be necessary to later split up $t=|S(I)|+|R|$. We also thus have %$|I|=2(t-d-1) \geq 2(q-1)(d+1)$.

%The same proof as Lemma 4.3(i) shows that the number of vertices in $S(I) %\cup R$ not in any $A_v$ is at most $2\Delta(t-d-1)-tn$.

%\textcolor{red}{Question: Why can we find an $(r-d-1)$-IT in $G$ omitting %exactly one vertex in each component of $I \cup R$? I understand that we can %find a partial IT by omitting exactly one vertex in each component of $S(I)$ %and then extending this to $G-R$ by Lemma 3.8, but I don't see why we can %find the vertices in $R$ to make this a full $(r-d-1)$-IT.}

%Assuming we can find an $(r-d-1)$-IT by the method in the question above, we %can now do something like the proof of Theorem 1.3. We assume that $G$ does %not have an $(r-d)$-IT. By Lemma 3.8, it suffices to show that we can find a %full IT $T_0$ in some component $J$ of $I$ such that $T_0 \subseteq \cup_{v %\in I} A_v$ This gives a contradiction as $T_0$ is independent of $I-(I \cap %J)$, so we can find an IT missing exactly one vertex in each of the other %components and then extend this to $G$ using Lemma 1.8 which gives a full $(r-%d)$-IT in $G$.

\section*{Acknowledgments}
The second author would like to thank the third author for generously providing mentorship and support over the past two years throughout this project.

\appendix

\section{Proof of Lemma \ref{algorithm terminates}}\label{appendix proof of lemma}
We prove Lemma \ref{algorithm terminates} in this section using similar arguments to those used in the proof of \cite[Theorem 2.2]{HS06}.

\begin{proof}[Proof of Lemma \ref{algorithm terminates}]
We will first show that Algorithm \ref{alg} maintains a feasible pair. Let $w, I, T$, and $T'$ be as defined in Step 2 or Step 3 of the algorithm. We will show that $(I \cup \{w\} \cup N(w,T'), T')$ is feasible. For convenience, let $I'=I \cup \{w\} \cup N(w,T')$. Also, we say that a class $V_i$ is \emph{active} in $I$ if $V_i \in S(I)$, i.e. $V_i \cap I \neq \emptyset$, and we refer to $S(I)$ as the set of active classes of $I$. It follows from the algorithm that $T \cap (\cup_{V_i \in S(I)} V_i)=T' \cap (\cup_{V_i \in S(I)} V_i)$, so that $T$ and $T'$ agree on active classes of $I$.                            

\textbf{Case 1:} We are at Step 2. 

Then (a) is satisfied since $T' \in \T$. 

For (b) suppose that $v \in T'$ and $S(\{v\}) \in S(I')=S(I) \cup S(N(w,T'))$. If $S(\{v\}) \in S(I)$ then $v \in T$ as $T$ and $T'$ agree on active classes of $I$. Then since $v \in T$ and $S(\{v\}) \in S(I)$, we have $v \in I \subseteq I'$ because $(I,T)$ satisfies (b). Now assume that $S(\{v\}) \in S(N(w,T'))$. Then we must have $v \in N(v,T') \subseteq I'$ as $T'$ is an IT.

We will now verify (c). By definition of $\T$, we have $W'=I' \setminus T'=W \cup \{w\}$ and $\{v\} \cup N(v,T')=\{v\} \cup N(v,T)$ for all $v \in W$. These stars $\{v\} \cup N(v,T')$ for $v \in W$ are all disjoint and of order at least two as $(I,T)$ is feasible. Also, $(\{w\} \cup N(w,T')) \cap (\{v\} \cup N(v,T))=\emptyset$ for all $v \in W$. Hence it only remains to be shown that $\{w\} \cup N(w,T')$ is a star with at least one leaf, i.e. $N(w,T') \neq \emptyset$. Suppose for the sake of contradiction that $N(w,T')=\emptyset$. Then $T'$ contains a vertex $u$ in the class of $w$ as if not, $T' \cup \{w\}$ is an IT in $G$ of size $|T'|+1$. But, since the class containing $w$ and $u$ intersects $I$, we have since $T$ and $T'$ agree on classes active in $I$, $u \in T$. Then $u \in I$ as well by (b) applied to $(I,T)$. Also, $u$ has exactly one neighbor in $I$ by (c) applied to $(I,T)$ as $u \not\in W$. Let $v_0 \in W$ be its neighbor. Then $T''=T' \cup \{w\} \setminus \{u\}$ is a partial IT on $S(T_0)$ such that $T'' \cap W=\emptyset$. We claim that the existence of $T''$ contradicts $(I,T)$ satisfying (e), so $N(w,T') \neq \emptyset$ and $(I',T')$ satisfies (c). We have $|N(v_0,T'')|=|N(v_0,T')\setminus\{u\}|=|N(v_0,T')|-1=|N(v_0,T)|-1$ and $N(v,T'')=N(v,T')=N(v,T)$ for $v \in W \setminus \{v_0\}$, giving us the desired contradiction.

Observe that $w \in V_I$ and $N(w,T') \cap V_I = \emptyset$. Hence $\G_{I'}$ is $\G_I$ with classes containing vertices in $N(w,T')$ added as leaves and is thus also a forest, so $(I',T')$ satisfies (d).

To show that $(I',T')$ satisfies (e), suppose for the sake of contradiction that there exists $v_0 \in W \cup \{w\}$ and $T'' \in \T'$ such that $T'' \cap W' = \emptyset$, $|N(v_0,T'')|<|N(v_0,T')|$, and $N(v,T'')=N(v,T')$ for all $v \in W' \setminus \{v_0\}$. Suppose that $v_0 \in W$. Then, since $N(v,T')=N(v,T)$ for all $v \in W$, $(I,T)$ would not satisfy (e), a contradiction. Hence we may assume $v_0=w$. But then $T'' \in \T$ and $\deg(w,T'')<\deg(w,T')$, which contradicts our choice of $T'$.

\textbf{Case 2:} We are at Step 3. In this case, we have $S(\{w\}) \not\in S(T)$. 

We have that (a) is satisfied since $T' \in \T$. 

Suppose that $v \in T'$ and $S(\{v\}) \in S(I')=S(I) \cup S(\{w\}) \cup S(N(w,T'))$. Since $w \not\in T'$ we have $S(\{v\}) \in S(I) \cup S(N(w,T'))$, and the same argument used in Case 1 above shows that $(I',T')$ satisfies (b).

We will now verify (c). By definition of $\T$, we have $W'=I' \setminus T'=W \cup \{w\}$ and $\{v\} \cup N(v,T')=\{v\} \cup N(v,T)$ for all $v \in W$. These stars $\{v\} \cup N(v,T')$ for $v \in W$ are all disjoint and of order at least two as $(I,T)$ is feasible. Also, $(\{w\} \cup N(w,T')) \cap (\{v\} \cup N(v,T))=\emptyset$ for all $v \in W$. Hence it only remains to be shown that $N(w,T') \neq \emptyset$. Suppose for the sake of contradiction that $N(w,T')=\emptyset$. But since $w \not\in S(T)$, this would imply that $\{w\} \cup T'$ is an IT in $G$ of size larger than $T$, a contradiction.

We have that $\G_{I'}$ is $\G_I$ with possibly a new vertex added for the class containing $w$ (if $S(\{w\}) \not\in S(I)$), and leaves corresponding to the classes containing vertices in $N(w,T')$ added. So, since $\G_I$ is a forest, so is $\G_{I'}$, and $(I',T')$ hence satisfies (d).

The same argument made in Case 1 above shows that $(I',T')$ satisfies (e).

We have thus showed that the algorithm maintains a feasible pair throughout. Since the size of $I$ increases throughout the algorithm, the algorithm must eventually terminate and output a feasible pair.
%conditions 1, 2, 4, and 5 follow immediately. We can see that $G[I]$ is a forest with a single component: the star $\{w\} \cup N(w,T')$. Hence to show that $(I', T')$ satisfies condition 3 it only remains to show that $N(w, T') \neq \emptyset$. Since if $N(w,T')=\emptyset$, $w \cup \{T'\}$ is an IT in $G$ with size larger than $|T'|$, we must have $N(w,T') \neq \emptyset$, so $(I', T')$ is feasible.
\end{proof}
\end{document}